\documentclass[11pt]{article}
\usepackage{graphicx}
\usepackage{color}
\usepackage{amsmath}
\usepackage{amssymb}
\usepackage{amscd}
\usepackage{bbm}
\newcommand{\R}{\mathbb{R}}
\newcommand{\inr}[1]{\bigl< #1 \bigr>}

\newcommand{\E}{\mathbb{E}}

\newcommand{\eps}{\varepsilon}

\newtheorem{Theorem}{Theorem}[section]
\newtheorem{Lemma}[Theorem]{Lemma}

\newtheorem{Definition}[Theorem]{Definition}

\newtheorem{Corollary}[Theorem]{Corollary}
\newtheorem{Remark}[Theorem]{Remark}

\newtheorem{Assumption}[Theorem]{Assumption}
\newtheorem{Question}[Theorem]{Question}
\numberwithin{equation}{section} 
\newcommand{\ls}{\leqslant}
\newcommand{\gr}{\geqslant}
\def \proof {\noindent {\bf Proof.}\ \ }

\def \endproof
{{\mbox{}\nolinebreak\hfill\rule{2mm}{2mm}\par\medbreak}}

\begin{document}
\title{{On generic chaining and the smallest singular value of random matrices with heavy tails}}
\author{Shahar Mendelson ${}^{1,3,4}$  \and Grigoris Paouris${}^{2,5}$}

\footnotetext[1]{Department of Mathematics, Technion, I.I.T, Haifa
32000, Israel.}
\footnotetext[2]{Department of Mathematics, Texas A$\&$M University, College Station, TX 77843-3368, U.S.A.}
\footnotetext[3]{Part of this research was supported by the Centre for Mathematics and its Applications,
The Australian National University, Canberra, ACT 0200,
Australia. Additional support was given by an Australian Research Council Discovery grant DP0559465, the European Community's Seventh Framework Programme (FP7/2007-2013) under ERC grant agreement 203134, and by the Israel Science Foundation grant 900/10.}
 \footnotetext[4] {Email:
shahar@tx.technion.ac.il}
\footnotetext[5] {Email:
grigoris@math.tamu.edu}

\maketitle

\begin{abstract}
We present a very general chaining method which allows one to control the supremum of the empirical process $\sup_{h \in H} |N^{-1}\sum_{i=1}^N h^2(X_i)-\E h^2|$ in rather general situations. We use this method to establish two main results. First, a quantitative (non asymptotic) version of the classical Bai-Yin Theorem on the singular values of a random matrix with i.i.d entries that have heavy tails, and second, a sharp estimate on the quadratic empirical process when $H=\{\inr{t,\cdot} : t \in T\}$, $T \subset \R^n$ and $\mu$ is an isotropic, unconditional, log-concave measure.
\end{abstract}

\section{Introduction} \label{sec:intro}
The main goal of this article is to obtain a non-asymptotic version of the Bai-Yin Theorem \cite{BaiYin} on the largest and smallest singular values of certain random matrices. The Bai-Yin theorem asserts the following:
\begin{Theorem} \label{thm:Bai-Yin}
Let $A=A_{N,n}$ be an $N \times n$ random matrix with independent entries, distributed according to a random variable $\xi$, for which
$$
\E \xi =0, \ \ \E \xi^2 =1 \ \ \E \xi^4 < \infty.
$$
If $N,n \to \infty$ and the aspect ratio $n/N$ converges to $\beta \in (0,1]$, then
$$
\frac{1}{\sqrt{N}} s_{\min}(A) \to 1 -\sqrt{\beta}, \ \ \
\frac{1}{\sqrt{N}} s_{\max}(A) \to 1 +\sqrt{\beta},
$$
almost surely, where $s_{\max}$ and $s_{\min}$ denote the largest and smallest singular value of $A$.

Also, without the fourth moment assumption, $s_{\max}(A) /\sqrt{N}$ is almost surely unbounded.
\end{Theorem}
The main result of this article is a quantitative version of the Bai-Yin Theorem.
\vskip0.5cm
\noindent {\bf Quantitative Bai-Yin Theorem.} For every $q>4$ and $L>0$, there exist constants $c_1$, $c_2$, $c_3$ and $c_4$ that depend only on $q$ and $L$ for which the following holds. For every integer $n$, $\beta \in (0,1]$ and $N=n/\beta$, let $A=A_{N,n}=(\xi_{i,j})$ be an $N \times n$ random matrix with independent, symmetric entries, distributed according to a random variable $\xi$, satisfying
 $\E \xi^2 =1$ and $\E |\xi|^q \leq L$.
Then, for any $n \geq c_1$, with probability at least $1-c_2/(\beta n^{c_3})$,
$$
1-c_4 \sqrt{\beta} \leq \frac{1}{\sqrt{N}} s_{\min}(A) \leq
\frac{1}{\sqrt{N}} s_{\max}(A) \leq 1 +c_4\sqrt{\beta}.
$$

The proof of this result is based on the analysis of a more general scenario which has been studied extensively in recent years, in which the given matrix has independent rows, selected according to a reasonable measure on $\R^n$, rather than a matrix with i.i.d. entries; and unlike the classical random matrix theory approach, one is naturally interested in the non-asymptotic behavior of the largest and smallest singular values of $\Gamma=N^{-1/2}\sum_{i=1}^N \inr{X_i,\cdot}e_i$ as a function of $N$ and $n$. We refer the reader to the surveys \cite{Ver-survey,RV-ICM} and references therein for the history and recent developments in the non-asymptotic theory of random matrices.

We will focus on the following questions:
\begin{Question} \label{Qu:main}
Let $\mu$ be a symmetric measure on $\R^n$ and let $(X_i)_{i=1}^N$ be selected independently according to $\mu$.
\item{1.} Let $\Sigma_N=\frac{1}{N}\sum_{i=1}^N X_i \otimes X_i$ be the sample covariance matrix and $\Sigma = \E (X \otimes X)$. Given $\eps>0$, is it true that with high probability, if $N \geq c(\eps)n$ then $\|\Sigma_N - \Sigma\|_{2 \to 2} \leq \eps$?
\item{2.} If $X$ is an isotropic vector (that is, $\E\inr{X,x}^2=\|x\|^2_{\ell_2^n}$ for every $x \in \R^n$), are there ``canonical" high probability bounds on $s_{\max}(\Gamma)$ and $s_{\min}(\Gamma)$? For example, under what conditions on $\mu$ are $s_{\max}(\Gamma)$ and $s_{\min}(\Gamma)$ of the order of $1 \pm c\sqrt{n/N}$ -- like in the Bai-Yin Theorem?
\end{Question}

Observe that the two questions are very similar. For example, it is straightforward to verify that if $\mu$ is isotropic, then both parts can be resolved by estimating the supremum of the empirical process
\begin{equation} \label{eq:emp-process}
\sup_{t \in S^{n-1}} \left|\frac{1}{N} \sum_{i=1}^N \inr{X_i,t}^2 -\E \inr{X,t}^2 \right|.
\end{equation}
And, in view of the second part of Question \ref{Qu:main}, we will be especially interested in the case $N \sim n$, that is, while keeping the aspect ratio $n/N$ constant.

When studying measures on $\R^n$ in this context, it is natural to divide the assumptions into two types: one on the $\ell_p^n$ norm of $X$ and the other on moments of linear functionals $\inr{x,\cdot}$.

To formulate the moment assumption we will use here, recall that for $\alpha \geq 1$, the $\psi_\alpha$ Orlicz norm of random variable $Z$ is defined by
$$
\|Z\|_{\psi_\alpha} = \inf \left\{c>0 : \E \exp(|Z|^\alpha/c^\alpha) \leq 2 \right\},
$$
and there are obvious extensions for $0<\alpha<1$. It is standard to verify that for every $\alpha>0$, $\|Z\|_{\psi_\alpha}$ is equivalent to $\sup_{q \geq 1} \|Z\|_{L_q}/q^{1/\alpha}$.

\begin{Assumption} \label{ass:moment-assumption}
For $p,q \geq 2$, a symmetric measure $\mu$ satisfies a $p$-small diameter, $L_q$ moment assumption with constants $\kappa_1$ and $\kappa_2$, if a random vector $X$ distributed according to $\mu$ satisfies that
\begin{equation} \label{eq:basic-assumption}
\|X\|_{\ell_p^n} \leq \kappa_1n^{1/p} \ {\rm a.s.,} \ \ {\rm and} \ {\rm for \ every} \ x \in S^{n-1}, \ \|\inr{x,\cdot}\|_{L_q} \leq \kappa_2.
\end{equation}
$\mu$ satisfies a small diameter $\psi_\alpha$ moment assumption if the $\psi_{\alpha}$ norm replaces the $L_q$ one in \eqref{eq:basic-assumption}.
\end{Assumption}

One should note that with very few exceptions, both parts of Assumption \ref{ass:moment-assumption} are needed if one wishes to address Question \ref{Qu:main}.

The $p$-small diameter component, i.e. that $\|X\|_{\ell_p^n} \leq \kappa_1n^{1/p}$ almost surely, is rather standard. Although it does not hold as stated even for a vector with i.i.d. gaussian entries,
one may assume it without loss of generality unless $N$ is much larger than $n$. Indeed, in typical situations $Pr(\|X\|_{\ell_p^n} \geq tn^{1/p})$ decays very quickly both in $t$ and in $n$. Therefore, $\max_{i \leq N} \|X_i\|_{\ell_p^n}/n^{1/p}$ is bounded with very high probability, unless $N$ is considerably larger than $n$ (see Section \ref{sec:pre} for more details). Hence, if $N \sim n$, which is the range we shall be interested in, a conditioning argument allows one to make the $p$-small diameter assumption.

Question \ref{Qu:main} has been studied under the $2$-small diameter assumption. In \cite{Rud99}, Rudelson showed that if $\|X\|_{\ell_2^n} \leq \kappa_1 \sqrt{n}$ almost surely then for every $N \geq c_1n\log n$, with probability at least $0.99$,
\begin{equation} \label{eq:Rudelson}
1-c_2\sqrt{\frac{n \log n}{N}} \leq s_{\min}(\Gamma) \leq s_{\max}(\Gamma) \leq 1+c_2\sqrt{\frac{n \log n}{N}},
\end{equation}
and $c_1, c_2$ are constants that depend only on $\kappa_1$.

It is straightforward to verify that this bound is optimal by considering the uniform measure on the set of coordinate vectors $\{\sqrt{n}e_1,...,\sqrt{n}e_n\}$, which results in the coupon-collector problem. Thus, given $\eps>0$, one requires at least $c(\eps)n\log n$ random points to ensure that the sample covariance matrix $\eps$-approximates the true covariance. Of course, \cite{Rud99} does not lead to a nontrivial estimate in the second part of Question \ref{Qu:main}, i.e. if the aspect ratio $n/N \to \beta \in (0,1]$ and $n \to \infty$, and in particular, \eqref{eq:Rudelson} can not yield a Bai-Yin type of bound. Any hope of getting the desired bounds in Question \eqref{Qu:main} requires additional assumptions on $X$.

Turning to the moments component of Assumption \ref{ass:moment-assumption}, note that a bound on the $L_q$ moments of linear functionals means that $\|\inr{x,\cdot}\|_{L_q} \lesssim \|x\|_{\ell_2^n}$, and if, in addition, $X$ is isotropic, the norms are equivalent. Moreover, in a similar fashion, a $\psi_\alpha$ assumption combined with isotropicity implies that the $\psi_\alpha$ and $\ell_2^n$ norms are equivalent.

Consider a situation when one only assumes such a moment condition. It is standard to verify that under a $\psi_2$ assumption, in which linear functionals exhibit a $\kappa_2$-subgaussian tail behavior (i.e., $Pr(|\inr{X,x}| \geq t\kappa_2\|x\|_{\ell_2^n}) \leq 2\exp(-t^2/2)$), then with probability at least $1-2\exp(-c_3n)$,
$$
s_{\min}(\Gamma),s_{\max}(\Gamma) \in [1-c_4\sqrt{n/N},1+c_4\sqrt{n/N}]
$$

Indeed, a Bernstein type inequality shows that for each $x \in S^{n-1}$ and $0<t<1/\kappa_2$,
$Pr(|N^{-1}\sum_{i=1}^N \inr{X_i,x}^2-\E \inr{X,x}^2| \geq t ) \leq 2\exp(-c_5Nt^2)$. And, if one is to obtain an estimate on the empirical process \eqref{eq:emp-process}, one has to control a $1/2$ net on the sphere, which is of cardinality $\sim \exp(c_6 n)$. The tradeoff between the complexity of the indexing set and the concentration at hand shows that with the desired probability,
$\sup_{t \in S^{n-1}} |N^{-1}\sum_{i=1}^N \inr{X_i,t}^2 -\E \inr{X,x}^2| \lesssim \sqrt{n/N}$.

Unfortunately, when one has a weaker moment estimate than a $\psi_2$ one, the situation becomes considerably more difficult. The complexity of the set one has to control remains the same, but the individual concentration deteriorates, because $N^{-1}\sum \inr{X_i,x}^2$ does not exhibit a strong enough concentration around its mean to balance the concentration-complexity tradeoff at the level of $\sqrt{n/N}$. Therefore, with a weaker moment assumption than a $\psi_2$ one, a combination of individual tail bounds and a ``global" assumption, like the small diameter information, is required in both parts of Question \ref{Qu:main}.

One situation in which the process \eqref{eq:emp-process} has been studied extensively in the last 15 years is a small diameter, $\psi_1$ moment assumption. The motivation for considering this situation comes from Asymptotic Geometric Analysis and the theory of log-concave measures, which are measures that have a symmetric, log-concave density. They fit the framework at hand nicely, because an isotropic, log concave vector $X$ satisfies that $\|X\|_{\ell_p^n} \leq c_1n^{1/p}$ with probability at least $1-2\exp(-c_2n^{1/p})$. Indeed, the case $p=2$ was proved in \cite{Pao}, while for $p>2$ the result was recently established by Lata{\l}a in \cite{Lat-ord-stat}. Moreover, linear functionals exhibit a $\psi_1$ behavior (see, e.g. \cite{Gia-survey} for a survey on log-concavity).

Partial results in the isotropic, log-concave case have been obtain by Bourgain \cite{Bour}, yielding an estimate on the covariance operator for $N =c(\eps)n\log^3 n$, which was improved by Rudelson \cite{Rud99} to $N=c(\eps) n \log^2n$. Subsequent improvements were $N =c(\eps)n\log n$ for unconditional convex bodies in \cite{GiaHaTs} and for general log-concave measures in \cite{Pao}. Finally, the optimal estimate of $N=c(\eps)n$ was obtained for an unconditional, log-concave measures by Aubrun \cite{Aub}, and for an arbitrary log-concave measure in Adamczak et al. \cite{ALPT1,ALPT2}, where the following result was proved:
\begin{Theorem} \label{thm:ALPT}
There exist absolute constants $c_1$ and  $c_2$ for which the following holds. If $\mu$ is an isotropic, log-concave measure, then with probability at least $1-\exp(-c_1\sqrt{n})$,
$$
\sup_{t \in S^{n-1}} \left| \frac{1}{N}\sum_{i=1}^N \inr{X_i,t}^2 - 1 \right| \leq c_2 \sqrt{\frac{n}{N}}.
$$
\end{Theorem}

Naturally, Question \ref{Qu:main} becomes even harder when one assumes that linear functionals have heavy tails, because sums of independent random variable exhibit very limited concentration -- far below the level required for the proof of Theorem \ref{thm:ALPT}. Recently, Vershynin \cite{Ver-heavy} proved the following remarkable fact:
\begin{Theorem} \label{thm:vershynin}
For every $q>4$, $\delta>0$ and constants $\kappa_1$ and $\kappa_2$, there exist constants $c_1$ and $c_2$ that depend on $q$, $\delta$ and $\kappa_1, \kappa_2$ for which the following holds.

If $\mu$ satisfies a $2$-small diameter, $L_q$ moment assumption with constants $\kappa_1$ and $\kappa_2$, then for every $\delta>0$, with probability at least $1-\delta$,
$$
\|\Sigma_N - \Sigma\|_{2 \to 2} \leq  c_1(\log \log n)^2 \left(\frac{n}{N}\right)^{1/2-2/q}.
$$
In particular, if $\mu$ is isotropic then
$$
1-c_2 \left(\frac{n}{N}\right)^{1/2-2/q} (\log \log n)^2\leq s_{\min}(\Gamma) \leq s_{\max}(\Gamma) \leq 1+c_2\left(\frac{n}{N}\right)^{1/2-2/q}(\log \log n)^2.
$$
\end{Theorem}

\noindent Moreover, very recently Strivastava and Vershynin \cite{SV}, obtained the following result:
\begin{Theorem}\label{thm:SV}
For every $\eta>0$, $\eps>0$ and $\kappa>0$ there exists constants $c_{1}, c_{2}$ and $c_{3}=\frac{\eta}{2\eta+2}$ for which the following holds. Let $\mu$ be an isotropic measure, satisfying that for every projection $P$ in $\R^n$,
$$
Pr\{ \|PX\|_{2}^{2} >t \} \leq \frac{\kappa} {t^{1+\eta}} , \ {\rm for} \ t\geq \kappa \ {\rm rank}(P). \leqno(\ast)
$$
If $(X_{i})_{i=1}^N$ are independent random vectors distributed according to $\mu$
then for every $N\gr c_{1}n$,
$$
\mathbb E \| \Sigma_{N} - I_{d}\| \ls \eps.
$$
Moreover, only under a $q$-moment assumption,
$$
1-c_2 \left(\frac{n}{N}\right)^{c_{3}} \leq \mathbb E s_{\min}(\Gamma)
$$
\end{Theorem}

It should be noted that the boundedness assumption in Theorem \ref{thm:SV} is satisfied by a vector with independent components $X=(\xi_i)_{i=1}^n$, if $\xi \in L_q$ for $q>4$, and thus both parts may be used in the i.i.d situation. However, for any $\eta>0$, $c_{3}<\frac{1}{2}$ ($1/2$ being the power in the Bai-Yin Theorem).

Our main result gives a version of Theorem \ref{thm:vershynin} for an unconditional measure with ``heavy tails".

\noindent{\bf Theorem A.} Let $\mu$ be an unconditional measure that satisfies the $p$-small diameter, $L_q$ moment assumption with constants $\kappa_1$ and $\kappa_2$ for some $p>2$.
\begin{description}
\item{1.} For every $q>4$ and $\delta < 1/2-1/2(p-1)$, there exist constants $c_0$, $c_1$ and $c_2$ that depend on $q$, $p$, $\kappa_1$, $\kappa_2$ and $\delta$, such that, for every $n \leq N \leq \exp(c_0n^\delta)$, with probability at least $1-\exp(-c_1n^\delta)$,
    $$
\sup_{t \in B_2^n} |N^{-1}\sum_{i=1}^N \inr{X_i,t}^2 - \E \inr{X,t}^2| \leq c_2\left(\frac{n}{N}\right)^{1/2}.
    $$
\item{2.} For every $2<q \leq 4$, if $p>(1-2/q)^{-1}$ and $\delta < 1/2-1/2(p-1)$,
there exist constants $c_3$ and $c_4$ that depend on $q$, $p$, $\delta$, $\kappa_1$ and $\kappa_2$, such that, for every $n \leq N \leq \exp(c_0n^\delta)$, with probability at least $1-\exp(-c_3n^\delta)$,
    $$
\sup_{t \in B_2^n} |N^{-1}\sum_{i=1}^N \inr{X_i,t}^2 - \E \inr{X,t}^2| \leq c_4\left(\frac{n}{N}\right)^{1-2/q}\log(N/n).
    $$
\end{description}
In both cases, for every $\eps>0$, with probability at least $1-2\exp(-cn^\delta)$,
$\|\Sigma_N - \Sigma\|_{2 \to 2} \leq  \eps$ provided that $N \gtrsim_{q,p,\delta,\kappa_1,\kappa_2} n$. Moreover, if $\mu$ is isotropic and $q >4$, then
    $$
1-c_2\left(\frac{n}{N}\right)^{1/2} \leq s_{\min}(\Gamma) \leq s_{\max}(\Gamma) \leq 1+c_2\left(\frac{n}{N}\right)^{1/2},
    $$
and if $2<q \leq 4$ then
    $$
1-c_4\left(\frac{n}{N}\right)^{1-2/q}\log(N/n) \leq s_{\min}(\Gamma) \leq s_{\max}(\Gamma) \leq 1+c_4\left(\frac{n}{N}\right)^{1-2/q}\log(N/n).
    $$
\vskip0.5cm

Our quantitative version of the Bai-Yin Theorem follows from Theorem A, because of the straightforward observation that if $\xi \in L_q$ for $q>4$ and is symmetric, then $X=(\xi_i)_{i=1}^n$ is unconditional, and there is some $p>2$ for which $\max_{i \leq N} \|X\|_{\ell_p^n} \lesssim n^{1/p}$ with high enough probability. Thus, conditioning $\mu$ to the unconditional body $cn^{1/p}B_p^n$ yield the desired result.

The approach we take in the proof of Theorem A is very different from all the previous results mentioned above, as those rely heavily on the fact that the empirical process \eqref{eq:emp-process} is indexed by the sphere or by the Euclidean ball, and that the underlying class of functions consists of linear functionals. At the heart of the arguments are either the classical trace method \cite{Aub}, a non-commutative Khintchine inequality \cite{Rud99} or sharp estimates on $\max_{|I|=k}\|\sum_{i \in I} X_i\|_{\ell_2^n}$ \cite{Bour,ALPT1,Ver-heavy}. As such, all these proofs are ``Euclidean" in nature and can not lead to bounds on the empirical process
\begin{equation} \label{eq:emp-general}
\sup_{h \in H} \left|\frac{1}{N}\sum_{i=1}^N h^2(X_i)-\E h^2\right|
\end{equation}
 for an arbitrary class of functions $H$ -- not even for $H_T=\{\inr{t,\cdot} : t \in T\}$ when $T$ is not the sphere or close to the sphere in some sense.

One should note that process \eqref{eq:emp-general} is an interesting object in its own right. For example, it has a key role in analyzing the uniform central limit Theorem \cite{Dud:book}; and, when indexed by $H_T$ for $T \subset \R^n$, it appear naturally in Asymptotic Geometric Analysis, for example, when proving embedding results or ``low-$M^*$" estimates for various matrix ensembles (see \cite{Men-psi-1} for a more detailed discussion). Thus, understanding what governs \eqref{eq:emp-general}, and in particular, going beyond the case $H_{B_2^n}$ is rather important.

The proof of Theorem A does just that, since it is based on a bound on \eqref{eq:emp-general} in terms of a certain notion of ``complexity" of the class $H$. It is not tailored to the case $H_{B_2^n}$, nor does it relay on the fact that the indexing class consists of linear functionals. Rather, the proof is based on a chaining scheme which is much more general than the applications that will be presented here.

The second application we chose to present as an illustration of the potential this empirical processes based method has, is the following.

Let $y_1,...y_n$ be independent, standard exponential random variables (i.e., with density $\sim \exp(-\sqrt{2}|t|)$, and for every $T \subset \R^n$ set
$$
E(T)=\E \sup_{t \in T} \sum_{i=1}^n t_i y_i, \ \ \ \ d_2(T)=\sup_{t \in T} \|t\|_{\ell_2^n}.
$$

\noindent{\bf Theorem B.} There exists absolute constants $c_1$, $c_2$ and $c_3$ for which the following holds. If $\mu$ is an isotropic, unconditional, log-concave measure on $\R^n$ and $T \subset \R^n$ is centrally symmetric, then for every $u \geq c_1$, with probability at least $1-2\exp(-c_2u^2)$,
\begin{equation} \label{eq:uncond-est}
\sup_{t \in T} \left|\frac{1}{N} \sum_{i=1}^N \inr{t,X_i}^2- \|t\|_{\ell_2^n}^2\right| \leq c_3 u^3\left(\frac{E(T)}{\sqrt{N}}+\frac{(E(T))^2}{N}\right).
\end{equation}
\vskip0.5cm

To put Theorem B in the right context, recall that a symmetric measure $\nu$ on $\R^n$ $(\kappa,L)$-weakly dominates a symmetric measure $\mu$ if for every $x \in \R^n$, and every $t>0$, $Pr_\mu(|\inr{x,\cdot}| \geq L t) \leq \kappa  Pr_\nu(|\inr{x,\cdot}| \geq t)$ \cite{KwaWoy}. For example, if $\mu$ is an isotropic $L$-subgaussian measure and $G=(g_1,...,g_n)$ is a standard gaussian vector in $\R^n$ then
$$
Pr_\mu(|\inr{x,\cdot}| \geq  Lt ) \leq 2\exp(-t^2/2\|x\|^2_{\ell_2^n}) = Pr_G(|\inr{x,\cdot}| \geq t),
$$
and thus $\mu$ is weakly dominated by $G$.

By the Majorizing Measures Theorem (see, e.g., \cite{Tal-book} and Section \ref{sec:pre}), it follows that if $\mu$ is $L$-subgaussian, there is a constant $c=c(L)$ satisfying that for every $T \subset \R^n$ and every integer $N$,
\begin{equation} \label{eq:linear-gaussian}
\E \sup_{t \in T} \inr{\sum_{i=1}^N X_i,t} \leq c \E \sup_{t \in T} \inr{\sum_{i=1}^N G_i,t} \equiv c\sqrt{N}G(T)
\end{equation}
where $(X_i)_{i=1}^N$ are independent copies of $X$, $(G_i)_{i=1}^N$ are independent copies of $G$ and $G(T)=\E \sup_{t \in T} \inr{G,t}$.

Moreover, the results of \cite{MPT,Men-psi-1} show that if $T$ is centrally symmetric and $\mu$ is isotropic and $L$-subgaussian, then
\begin{equation} \label{eq:quadratic-gaussian}
\E \sup_{t \in T} \left|\frac{1}{N} \sum_{i=1}^N \inr{t,X_i}^2- \|t\|_{\ell_2^n}^2\right| \lesssim_L \frac{G(T)}{\sqrt{N}}+\frac{(G(T))^2}{N}.
\end{equation}
Hence, the fact that an $L$-subgaussian measure is weakly dominated by a gaussian measure (with the same covariance structure) is exhibited by a strong domination in \eqref{eq:linear-gaussian} and in \eqref{eq:quadratic-gaussian}, that holds for every $T \subset \R^n$.

Just like subgaussian vectors, isotropic, unconditional log-concave vectors  have a natural weakly dominating measure. By the Bobkov-Nazarov Theorem \cite{BobNaz} they are $(\kappa,L)$-weakly dominated by the vector $Y=(y_1,...,y_n)$, and $\kappa$ and $L$ are absolute constants. In \cite{Lat-weak}, Lata{\l}a showed that as in \eqref{eq:linear-gaussian}, for every $T \subset \R^n$,
$\E \sup_{t \in T} \inr{\sum_{i=1}^N X_i,t} \lesssim \E \sup_{t \in T} \inr{\sum_{i=1}^N Y_i,t}$. Theorem B shows that the quadratic strong domination, analogous to \eqref{eq:quadratic-gaussian}, is also true in this case.

Theorem B has many standard applications, leading to embedding results of a similar nature to the Johnson-Lindenstrauss Lemma and to ``low $M^*$" estimates that hold for unconditional, log-concave ensembles. Deriving these and other outcomes from Theorem B is standard and will not be presented here. One should also note that a log-concave Chevet type inequality, i.e., upper estimates on the operator norm $\|\Gamma\|_{X \to Y}$ for finite dimensional normed spaces $X$ and $Y$ has recently been established in \cite{ALLPT}.

In the next section we will present several preliminary facts and definitions that will be used throughout this article. Then, in Section \ref{sec:decomp} we will show that if $V \subset \R^N$ can be decomposed in a certain way, the Bernoulli process indexed by $\{(v_i^2)_{i=1}^N : v \in V\}$ is well behaved. Section \ref{sec:cord-proj} is devoted to the observation that if $H$ is a class of functions, then under mild assumptions and with high probability, the random coordinate projection $P_\sigma H =\{(h(X_i))_{i=1}^N : h \in H\}$ can be decomposed in the sense of Section \ref{sec:decomp}. It turns out that the decomposition depends on the complexity of $H$ and on the decay of tails of functions in $H$. Finally, in Section \ref{sec:estimates} we will present examples in which the complexity of $H$ can be estimated, leading to the proofs of Theorem A (and consequently, the quantitative Bai-Yin Theorem) and of Theorem B.

\section{Preliminaries} \label{sec:pre}
Throughout, all absolute constants are positive numbers, denoted by $c,c_0,c_1,...$ and their value may change from line to line.
$\kappa_0,\kappa_1,...$ denote constants whose value will remain
unchanged. By $A \sim B$ we mean that there are absolute constants
$c$ and $C$ such that $cB \leq A \leq CB$, and by $A \lesssim B$
that $A \leq CB$. $A \sim_\gamma B$ (resp. $A \lesssim_\gamma B$) denotes that the constants depend only on $\gamma$.

For $1 \leq p \leq \infty$, $\ell_p^n$ is $\R^n$
endowed with the $\ell_p$ norm, which we denote by $\| \ \|_{\ell_p^n}$, and
$B_p^n$ is its unit ball. With a minor abuse of notation we write
$| \ |$ both for the cardinality of a set and for the absolute value. Finally, if $(a_n)$ is a sequence, let $(a_n^*)$ be a non-increasing rearrangement of $(|a_n|)$.

Next, let us turn to the complexity parameters that motivated our method of analysis -- Talagrand's $\gamma$-functionals.
\begin{Definition} \label{def:gamma-2} \cite{Tal-book}
For a metric space $(T,d)$, an {\it admissible sequence} of $T$ is a
collection of subsets of $T$, $\{T_s : s \geq 0\}$, such that for
every $s \geq 1$, $|T_s| \leq 2^{2^s}$ and $|T_0|=1$. For $\beta \geq 1$,
define the $\gamma_\beta$ functional by
$$
\gamma_\beta(T,d) =\inf \sup_{t \in T} \sum_{s=0}^\infty
2^{s/\beta}d(t,T_s),
$$
where the infimum is taken with respect to all admissible sequences
of $T$. For an admissible sequence $(T_s)_{s \geq 0}$ we denote by
$\pi_s t$ a nearest point to $t$ in $T_s$ with respect to the metric
$d$.
\end{Definition}
One should note that our chaining approach is based on a slightly less restrictive definition, giving one more freedom; for example, the cardinality of the sets will not necessarily be $2^{2^s}$, the metric may change with $s$, etc. (see Section \ref{sec:decomp}).

When
considered for a set $T \subset L_2$, $\gamma_2$ has close connections with
properties of the canonical gaussian process indexed by $T$, and we
refer the reader to \cite{Dud:book,Tal-book} for detailed
expositions on these connections. One can show that under mild measurability
assumptions, if $\{G_t: t \in T\}$ is a centered gaussian process indexed by a set $T$, then
$$
c_1 \gamma_2(T,d) \leq \E \sup_{t \in T} G_t \leq c_2 \gamma_2(T,d),
$$
where $c_1$ and  $c_2$ are  absolute constants and for every $s,t
\in T$, $d^2 (s,t) = \E|G_s-G_t|^2$.  The upper bound is due to
Fernique \cite{F} and the lower bound is Talagrand's Majorizing
Measures Theorem \cite{Tal87}. Note that if $T \subset \R^n$,
$(g_i)_{i=1}^n$ are standard, independent gaussians and $G_t =
\sum_{i=1}^n g_i t_i$ then $d(s, t) = \|s-t\|_{\ell_2^n}$, and therefore
\begin{equation}
  \label{maj_meas}
  c_1 \gamma_2(T,\| \cdot \|_{\ell_2^n}) \leq  \E \sup_{t \in T} \sum_{i=1}^n
  g_it_i  \leq c_2 \gamma_2(T,\|\cdot \|_{\ell_2^n}).
\end{equation}

A part of our discussion (Theorem B) will be devoted to isotropic, log-concave measures on $\R^n$.

\begin{Definition}
A symmetric probability measure $\mu$ on $\R^n$ is isotropic
if for every $y \in \R^n$, $\int |\inr{x,y}|^2 d\mu(x) = \|y\|_{\ell_2^n}^2$.

The measure $\mu$ is log-concave if for every $0<\lambda<1$ and
every nonempty Borel measurable sets $A,B \subset \R^n$,
$\mu(\lambda A+(1-\lambda)B) \geq \mu(A)^\lambda\mu(B)^{1-\lambda}$.
\end{Definition}
A typical example of a log-concave measure on $\R^n$ is the volume
measure of a convex body in $\R^n$, a fact that follows from the
Brunn-Minkowski inequality (see, e.g. \cite{Pis:book}). Moreover,
Borell's inequality \cite{Bor,MilSch} implies that there is an
absolute constant $c$ such that if $\mu$ is an isotropic,
log-concave measure on $\R^n$, then for every $x \in \R^n$,
$\|\inr{x,\cdot}\|_{\psi_1} \leq c\|\inr{x,\cdot}\|_{L_2} =
c\|x\|_{\ell_2^n}$.

As mentioned in the introduction, if $X$ is distributed according to an isotropic, log-concave measure on $\R^n$ then $\|X\|_{\ell_p^n}$ decays quickly at scales that are larger than $n^{1/p}$. Thus, by conditioning, the main result in \cite{Pao} shows that a $2$-small diameter assumption can be made without loss of generality as long as $N \leq \exp(c\sqrt{n})$, and Lata{\l}a \cite{Lat-ord-stat} proved the analogous result for $p>2$, as long as $N \leq \exp(cn^{1/p})$.


\section{Decomposition of sets} \label{sec:decomp}
We begin with a description of the modified chaining procedure. Let $(\eta_s)_{s \geq 0}$ be an increasing sequence which satisfies that for every $s \geq 0$, $2^{\eta_s} \cdot 2^{\eta_{s+1}} \leq 10 \cdot 2^{\eta_{s+2}}$ and for $s \geq 1$, $1.1 \leq \eta_{s+1}/\eta_s \leq 10$ (where 1.1 can be replaced by $1+\eps$ and 10 can be any suitably large constant). For example, $\eta_0=0$ and $\eta_s=2^s$ for $s \geq 1$ is the usual choice of a sequence that has been used in the definition of Talagrand's $\gamma$ functionals. An admissible sequence of $V \subset \R^N$ relative to $(\eta_s)_{s \geq 0}$ is a collection of subsets $V_s \subset V$ for which $|V_s| \leq 2^{\eta_s}$. For every $s$ let $\pi_s:V \to V_s$, which usually will be a nearest point map relative to some distance. We will denote $\pi_s v - \pi_{s-1} v$ by $\Delta_s v$, and sometimes write $\Delta_0 v$ for  $\pi_0 v$. Finally, $\Delta_s V$ is the set $\{\Delta_s v : v \in V\}$.

Let $\phi$ be an increasing function which will be chosen according to additional information one will have on the given class. Examples that one should have in mind are $\phi_\beta(x) \sim_\beta \sqrt{x}\log^{1/\beta}(eN/x)$, resulting from a bound on the $\psi_\beta$ diameter of $H$, or $\phi_{q,\eps} \sim_{q,\eps} N^{(1+\eps)/q} x^{1/2-(1+\eps)/q}$ for $q>2$ and $\eps$ in the right range, arising from an $L_q$ moment assumption.

Assume that $V \subset \R^N$ is endowed with a family of functionals $\theta_s$ and a semi-norm $\| \ \|$ (which, in our applications, will either arise from the $L_q$ norm or from the $\psi_\beta$ norm), and set $d = \sup_{v \in V} \|v\|$.
\begin{Definition} \label{def:decomposition}
$V \subset \R^N$ admits a decomposition with constants $\alpha$ and $\gamma$ if it has an admissible sequence $(V_s)_{s \geq 0}$ relative to $(\eta_s)_{s \geq 0}$ for which the following holds.
\begin{description}
\item{1.} $\sup_{v \in V} \left(\theta_0(\pi_0 v)+\sum_{s > 0} \theta_s(\Delta_s v) \right) \leq \gamma$.
\item{2.} For every $v \in V$ and every $I \subset \{1,...,N\}$,
$$
\left(\sum_{i \in I} v_i^2 \right)^{1/2} \leq \alpha \left(\gamma+d \phi(|I|)\right).
$$
\item{3.} If  $\eta_s \leq N$ then for every $v \in V$ and every $I\subset \{1,...,N\}$
$$
\left(\sum_{i \in I} (\Delta_s v)^2_i \right)^{1/2} \leq \alpha\left(\theta_s(\Delta_s v)+\|\Delta_s v\| \phi(|I|)\right),
$$
and if $\eta_s \geq N$ then for every $v \in V$ and every $I \subset \{1,...,N\}$,
$$
\left(\sum_{i \in I} (\Delta_s v)^2_i \right)^{1/2} \leq \alpha \theta_s(\Delta_s v).
$$
\end{description}
\end{Definition}

Although this definition seems artificial at first glance, we will show that it captures the geometry of a typical coordinate projection $P_\sigma H =\{ (h(X_i))_{i=1}^N : h \in H\}$.

The main observation of this section is that one can use this type of decomposition to bound the supremum of the Bernoulli process indexed by $V^2=\{(v_i^2)_{i=1}^N : v \in V\}$. Hence, if $V=P_\sigma H$, then a standard symmetrization argument leads to the desired bound on $\sup_{h \in H} |N^{-1}\sum_{i=1}^N h^2(X_i) -\E h^2|$ (see section \ref{sec:proofsA-B}).

To formulate the estimate on the Bernoulli process, set
$$
\Phi=\left(\sum_{i=1}^N \frac{\phi^4(i)}{i^2}\right)^{1/2},  \ \ \
\Phi_s=\left(\sum_{i=1}^{N-\eta_s} \frac{\phi^2(\eta_s+i)}{\eta_s+i}\cdot \frac{\phi^2(i)}{i}\right)^{1/2}
$$
for $\eta_s \leq N$, put
$$
A_1=\sup_{v \in V} \sum_{\{s>0: \eta_s \leq N\}} \phi(\eta_s)\|\Delta_s v\|, \ \ A_2 = \sup_{v \in V} \sum_{\{s>0: \eta_s \leq N\}} \phi^2(\eta_s)\|\Delta_s v\|,
$$
and let
$$
A_\Phi=\sup_{v \in V} \sum_{\{s: \eta_s \leq N\}} \Phi_s \eta_s^{1/2}\|\Delta_s v\|.
$$
For $2 < q \leq 4$ and $0 \leq \eps<(q/2)-1$, let
$$
 B_{q,\eps}= \sup_{v \in V} \sum_{\{s: \eta_s \leq N\}} \eta_s^{1-2(1+\eps)/q}\|\Delta_s v\|.
$$
As will become clearer, the most important of the $B_{q,\eps}$ parameters is
$$
B_4 \equiv B_{4,0}=\sup_{v \in V} \sum_{\{s: \eta_s \leq N\}} \eta_s^{1/2}\|\Delta_s v\|,
$$
which, under the standard choice of $\eta_0=0$ and $\eta_s=2^s$ for $s \geq 1$, corresponds to $\gamma_2(V,\| \ \|)$.

\begin{Theorem} \label{thm:bernoulli-bounds}
There exist absolute constants $c_0$, $c_1$ and $c_2$ for which the following holds. If $V \subset \R^N$ has a decomposition as in Definition \ref{def:decomposition}, then for every $r \geq c_0$, with probability at least $1-2\exp(-c_1r^2\eta_0)$,
$$
\sup_{v \in V} \left|\sum_{i=1}^N \eps_i v_i^2 \right| \leq c_2 r \alpha^2
\left(\gamma(\gamma+d\phi(N)+A_1) + d\left(A_2+ A_\Phi \right)\right)
$$
\end{Theorem}

Before presenting the proof, let us consider the two main examples which will interest us, namely, the families $\phi_\beta=\sqrt{x}\log^{1/\beta}(eN/x)$ for any $\beta>0$ and $\phi_{q,\eps}=\sqrt{x}(N/x)^{(1+\eps)/q}$ for any $q>2$ (and for $\eps$ selected appropriately).

In both cases $\phi(N) \sim \sqrt{N}$ and for any $\beta>0$, $\Phi \sim_\beta \sqrt{N}$. If $q>4$ and $0 \leq \eps \leq q/4-1$, $\Phi \leq (1-4(1+\eps)/q)^{-1/2}\sqrt{N}$, and since $\Phi_s \leq \Phi$, then for $\beta>0$ or $q>4$,
$$
A_\Phi \leq \Phi \sup_{v \in V} \sum_{\{s: \eta_s \leq N\}} \eta_s^{1/2}\|\Delta_s v\| \lesssim \sqrt{N}B_4,
$$
with the constant depending either on $\beta$ or on $q$ and $\eps$ as above.

On the other hand, if $2<q \leq 4$ and $0<\eps<q/2-1$ then
\begin{align*}
\Phi_s = & \left(\sum_{i=1}^{N-\eta_s} \frac{\phi^2(\eta_s+i)}{\eta_s+i}\cdot \frac{\phi^2(i)}{i}\right)^{1/2}=\left(\sum_{i=1}^{N-\eta_s} \left(\frac{N}{\eta_s+i}\right)^{2(1+\eps)/q} \cdot \left(\frac{N}{i}\right)^{2(1+\eps)/q} \right)^{1/2}
\\
\lesssim & \frac{N^{2(1+\eps)/q}}{1-2(1+\eps)/q} \cdot  \eta_s^{1/2-2(1+\eps)/q}.
\end{align*}
Therefore, in that range
\begin{equation*}
A_\Phi \lesssim \frac{N^{2(1+\eps)/q}}{1-2(1+\eps)/q} \cdot \sup_{v \in V} \sum_{\{s: \eta_s \leq N\}} \eta_s^{1-2(1+\eps)/q}\|\Delta_s v\| = \frac{N^{2(1+\eps)/q}}{1-2(1+\eps)/q}B_{q,\eps}.
\end{equation*}

Next, since $(\eta_s)_{s \geq 0}$ increases exponentially, then for $q>2$
\begin{equation} \label{eq:eta}
\sum_{\{s: \eta_s \leq N\}} \phi(\eta_s) \|\Delta_s v\| \leq 2d \sum_{\{s: \eta_s \leq N\}} \phi(\eta_s) \lesssim d \sqrt{N},
\end{equation}
and the constant in \eqref{eq:eta} depends on $\beta$ or on $q$ and $\eps$ respectively. In particular, if $2<q \leq 4$ and $0<\eps<q/2-1$, then
$$
\sum_{\{s: \eta_s \leq N\}} \phi(\eta_s) \lesssim \frac{\sqrt{N}}{1-2(1+\eps)/q}
$$

Finally, one has to control $\sum_{\{s: \eta_s \leq N\}} \phi^2(\eta_s)\|\Delta_s v\|$. Note that if $\beta>0$ or $q>4$, then
\begin{align*}
\sum_{\{s: \eta_s \leq N\}} \phi^2(\eta_s)\|\Delta_s v\| \leq & \left(\max_{\{s : \eta_s \leq N\}} \frac{\phi^2(\eta_s)}{\eta_s^{1/2}}\right) \cdot \sum_{\{s: \eta_s \leq N\}} \eta_s^{1/2} \|\Delta_s v\|
\\
\lesssim & \sqrt{N} \sum_{\{s: \eta_s \leq N\}} \eta_s^{1/2} \|\Delta_s v\| \sim \sqrt{N} B_4,
\end{align*}
and if $2<q \leq 4$ and $0<\eps<q/2-1$ then
$$
\sum_{\{s: \eta_s \leq N\}} \phi^2(\eta_s)\|\Delta_s v\| \leq N^{2(1+\eps)/q} \sum_{\{s: \eta_s \leq N\}} \eta_s^{1-2(1+\eps)/q}\|\Delta_s v\| = N^{2(1+\eps)/q} B_{q,\eps}.
$$
We thus arrive to a more compact formulation of Theorem \ref{thm:bernoulli-bounds} in the cases we will be interested in.
\begin{Corollary} \label{cor:Bernoulli-formulation}
For any $\beta>0$ or $q>4$, with probability at least $1-2\exp(-c_1r^2\eta_0)$,
$$
\sup_{v \in V} \left|\sum_{i=1}^N \eps_i v_i^2 \right| \lesssim r \alpha^2  \left(
\gamma^2+d \sqrt{N}(\gamma+B_4)\right),
$$
with a constant that depends on $\beta$ or on $q$ and $\eps$ respectively.

Also, if $2<q \leq 4$ and $0<\eps<q/2-1$, then with probability at least $1-2\exp(-c_1r^2\eta_0)$,
$$
\sup_{v \in V} \left|\sum_{i=1}^N \eps_i v_i^2 \right| \lesssim \frac{\alpha^2 r}{1-2(1+\eps)/q} \left(
\gamma^2+d \sqrt{N}\gamma +d N^{2(1+\eps)/q} B_{q,\eps}\right).
$$
\end{Corollary}

\noindent {\bf Proof of Theorem \ref{thm:bernoulli-bounds}.} For every $\Delta_s v $ let $i$ be the largest integer in $\{1,...,N\}$ for which $\theta_s(\Delta_s v) \geq \|\Delta_s v\| \phi(i)$. Throughout the proof we will assume that such an integer exists, and if it does not, the necessary modifications to the proof are obvious. Let $i_{s,v}=\max\{i,\eta_s\}$ and put $I_{s,v}$ to be the set of the largest $i_{s,v}$ coordinates of $|\Delta_s v|$. Let $\Delta_s^+v=P_{I_{s,v}} \Delta_s v$ and $\Delta_s^-v=P_{I^c_{s,v}} \Delta_s v$ be the projections of $\Delta_s v$ onto the set of coordinates $I_{s,v}$ and $I_{s,v}^c$ respectively. Also, let $j$ be the largest integer in $\{1,...,N\}$ for which $\gamma \geq d \phi(j)$. Thus, for every $v \in V$,
$\left(\sum_{i=1}^j (v_i^2)^* \right)^{1/2} \leq 2\alpha \gamma$ and for every $\ell \geq j$, $v_\ell^* \leq 2\alpha d \phi(\ell)/\sqrt{\ell}$. If $J$ is the set of the largest $j$ coordinates of $v \in V$, let $v^+=P_J v$ and $v^-=P_{J^c} v$.

Let $w \cdot v = \sum_{i=1}^N w_i v_i e_i$, and since
$$
v^2-(\pi_0 v)^2=\sum_{s > 0} (\pi_{s} v)^2 - (\pi_{s-1} v)^2 = \sum_{s>0}  (\Delta_s v) \cdot (\pi_s v + \pi_{s-1} v),
$$
one has to control increments of the form
$\sum_{i=1}^N \eps_i (\Delta_s v )_i (\pi_s v+\pi_{s-1}v)_i$.

Observe that if $\eta_s \geq N$ then with probability $1$,
\begin{align*}
& \sum_{i=1}^N \eps_i ((\Delta_s v) \cdot (\pi_s v +\pi_{s-1} v))_i \leq
\| (\Delta_s v) \cdot (\pi_s v +\pi_{s-1} v)\|_{\ell_1^N}
\\
\leq &
2\|\Delta_s v\|_{\ell_2^N} \sup_{v \in V} \|v\|_{\ell_2^N}
\leq 2\alpha^2 \theta_s(\Delta_s v) (\gamma+d \phi(N)).
\end{align*}
Next, if $\eta_s \leq N$ we will decompose the vectors one has to control according to the size of their coordinates, because, with probability $1-2\exp(-r^2/2)$,
\begin{align} \label{eq:basic-subgaussian-Bernoulli}
\nonumber  & \left|\sum_{i=1}^N \eps_i (\Delta_s v)_i  (\pi_s v + \pi_{s-1} v)_i \right| \leq \|(\Delta_s^+ v) \cdot (\pi_s v + \pi_{s-1}v) \|_{\ell_1^N}
\\
+ & r  \| (\Delta_s^- v) \cdot ((\pi_s v)^+ + (\pi_{s-1} v)^+) \|_{\ell_2^N}
+ r\|(\Delta_s^- v) \cdot ((\pi_s v)^- + (\pi_{s-1} v)^-) \|_{\ell_2^N}.
\end{align}
Consider the following two cases. If $i_{s,v}=\eta_s$ then
\begin{align*}
\|(\Delta^+_s v) \cdot  (\pi_s v + \pi_{s-1} v)\|_{\ell_1^N} \leq & \|\Delta_s^+ v\|_{\ell_2^N} \|P_{I_{s,v}} (\pi_s v + \pi_{s-1} v)\|_{\ell_2^N}
\\
\lesssim_{\alpha^2} &  \|\Delta_s v\| \phi(\eta_s) \left(\gamma+d \phi(\eta_s) \right).
\end{align*}
Moreover,
$$
\|\Delta_s^- v\|_{\ell_\infty^N} \leq \frac{\|\Delta^+_s v\|_{\ell_2^N}}{|I_{s,v}|^{1/2}} \leq 2 \alpha \|\Delta_s v\| \frac{\phi(\eta_s)}{\eta_s^{1/2}},
$$
and thus, for every $v \in V$,
\begin{equation*}
\eta_s^{1/2} \| (\Delta_s^- v) \cdot   w^+\|_{\ell_2^N} \leq \eta_s^{1/2} \| \Delta_s^- v\|_{\ell_\infty^N} \|w^+\|_{\ell_2^N}
\lesssim_{\alpha^2} \gamma \phi(\eta_s)\|\Delta_s v\|.
\end{equation*}
To estimate $\eta_s^{1/2}\|(\Delta_s^- v) \cdot w^-\|_{\ell_2^N}$, observe that since
$(\Delta^{-}_s v)_i^* \lesssim_\alpha \|\Delta_s v\| \phi(\eta_s+i)/\sqrt{\eta_s +i}$, $w_i^* \lesssim_\alpha d \phi(i)/\sqrt{i}$ and $\sum |a_i b_i| \leq \sum a_i^* b_i^*$, then
\begin{align*}
\eta_s^{1/2}\|(\Delta_s^- v) \cdot w^-\|_{\ell_2^N} & \lesssim_{\alpha^2}
\eta_s^{1/2} \|\Delta_sv\| d \left(\sum_{i=1}^{N-\eta_s} \frac{\phi^2(i)}{i} \cdot \frac{\phi^2(\eta_s+i)}{\eta_s +i}\right)^{1/2}
\\
& \lesssim_{\alpha^2}d
\eta_s^{1/2} \Phi_s \|\Delta_sv\|.
\end{align*}

Therefore, summing the three terms over $\{s>0 : \eta_s \leq N\}$,
\begin{align*}
& \sum_{\{s>0 : \eta_s \leq N\}} \|(\Delta_s^+ v) \cdot (\pi_s v + \pi_{s-1} v)\|_{\ell_1^N}
\\
\lesssim_{\alpha^2} & \gamma \sum_{\{s>0 : \eta_s \leq N\}} \phi(\eta_s)\|\Delta_s v\| + d \sum_{\{s>0 : \eta_s \leq N\}}
\phi^2(\eta_s) \|\Delta_s v\|,
\end{align*}
\begin{equation*}
\sum_{\{s>0 : \eta_s \leq N\}} \eta_s^{1/2} \|(\Delta_s^- v) \cdot ((\pi_s v)^+ + (\pi_{s-1} v)^+)\|_{\ell_2^N}
\lesssim_{\alpha^2} \gamma \sum_{\{s>0 : \eta_s \leq N\}} \phi(\eta_s)\|\Delta_s v\|,
\end{equation*}
and
\begin{equation*}
\sum_{\{s>0 : \eta_s \leq N\}} \eta_s^{1/2} \|(\Delta_s^- v) \cdot ((\pi_s v)^- + (\pi_{s-1} v)^-)\|_{\ell_2^N}
\lesssim_{\alpha^2} d  \sum_{\{s>0 : \eta_s \leq N\}} \eta_s^{1/2}\Phi_s \|\Delta_sv\|.
\end{equation*}

Next, if $i_{s,v} \not = \eta_s$ then $\|\Delta_s^+ v \|_{\ell_2^N} \leq 2\alpha \theta_s(\Delta_s v)$, and thus
$$
\|(\Delta_s^+ v) \cdot (\pi_s v + \pi_{s-1} v)\|_{\ell_1^N} \leq 2\alpha^2 \theta_s(\Delta_s v) (\gamma+d \phi(N)).
$$
Since $|I_{s,v}| \geq \eta_s$,
$$
\|\Delta_s^- v\|_{\ell_\infty^N} \leq 2\alpha  \theta_s(\Delta_s v) /|I_{s,v}|^{1/2} \leq 2\alpha  \frac{\theta_s(\Delta_s v)}{\eta_s^{1/2}},
$$
then splitting each $w \in V$ to $w^+ + w^-$ as above,
\begin{equation*}
\eta_s^{1/2} \|(\Delta_s^- v) \cdot w^+\|_{\ell_2^N}
\leq \eta_s^{1/2} \|\Delta_s^- v\|_{\ell_\infty^N} \|w^+\|_{\ell_2^N}
\lesssim_{\alpha^2} \gamma \theta_s(\Delta_s v),
\end{equation*}
and
\begin{equation*}
\eta_s^{1/2} \|(\Delta_s^- v) \cdot w^-\|_{\ell_2^N} \lesssim_{\alpha^2}
d \eta_s^{1/2} \Phi_s \|\Delta_sv\|.
\end{equation*}
Therefore,
\begin{align*}
(*) =& \sum_{\{s>0 : \eta_s \leq N\}} \|(\Delta_s^+ v) \cdot (\pi_s v + \pi_{s-1} v)\|_{\ell_1^N} + \eta_s^{1/2} \|(\Delta_s^- v)^ \cdot ((\pi_s v)^+ + (\pi_{s-1}v)^+)\|_{\ell_2^N}
\\
+ & \eta_s^{1/2} \|(\Delta_s^- v) \cdot ((\pi_s v)^- + (\pi_{s-1}v)^-)\|_{\ell_2^N}
\leq (3)+(4)+(5),
\end{align*}
where
\begin{equation*}
(3) \lesssim_{\alpha^2} (\gamma+d\phi(N)) \sum_{\{s>0 : \eta_s \leq N\}} \theta_s(\Delta_s v), \ \ \
(4) \lesssim_{\alpha^2} \gamma \sum_{\{s>0 : \eta_s \leq N\}} \theta_s(\Delta_s v),
\end{equation*}
and
\begin{equation*}
(5) \lesssim_{\alpha^2} d \sum_{\{s>0 : \eta_s \leq N\}} \eta_s^{1/2} \Phi_s \|\Delta_sv\|.
\end{equation*}

Recall that $|\Delta_s V|, |V_s| \leq 10 \cdot 2^{\eta_{s+1}}$ and that $\eta_{s+1} \leq 10 \eta_s$. Given $r \geq c_0$, then applying \eqref{eq:basic-subgaussian-Bernoulli} for $t_s =10 r\eta_s^{1/2}$ and summing over $\{s: \eta_s \leq N\}$, it follows that $\sup_{v \in V} \left|\sum_{i=1}^N \eps_i (v^2-(\pi_0 v)^2)_i \right|$ is bounded by the desired quantity with probability at least $1-2\exp(-c_1r^2\eta_0)$.

Finally, for $v \in V_0$, let $i$ be the largest integer in $\{1,...,N\}$ for which  $\theta_0(v) \geq \|v\| \phi(i)$, and set $I$ to be the set of the $i$-largest coordinates of $v$. Thus, $\sum_{i \in I} v_i^2 \leq 2\alpha^2 \theta_0(v) \leq 2\alpha^2 \gamma^2$, and for $\ell \geq i$, $v_\ell^* \leq \alpha \|v\| \phi(\ell)/\sqrt{\ell}$. Since $|V_0| \leq 2^{\eta_0}$, then with probability at least $1-2\exp(-c_2r^2\eta_0)$
$$
\left|\sum_{i=1}^N \eps_i v_i^2 \right| \lesssim_{\alpha^2}  \gamma^2+rd \Phi_0 \eta_0^{1/2}\|v\|,
$$
completing the proof.
\endproof

\section{Coordinate projections of Function classes} \label{sec:cord-proj}
The aim of this section is to show that under very mild assumptions, empirical processes have well behaved coordinate projections in the sense of Definition \ref{def:decomposition}. A first result in this direction was established in \cite{Men-psi-1}, in which the main observation, formulated in the language of Section \ref{sec:decomp}, was that if $\eta_0=0$ and $\eta_s = 2^s$ for $s \geq 1$, then for the choice of $\theta_{s}((h(X_i))_{i=1}^N)=2^{s/2}\|h\|_{\psi_2}$, $\alpha=\sqrt{u}$, $\|(h(X_i))_{i=1}^N\|=\|h\|_{\psi_1}$ and $\phi(x) \sim \sqrt{x} \log(eN/x)$, the set $V=\{(h(X_i))_{i=1}^N : h \in H\}$ has a good decomposition with high probability.
Hence, the Bernoulli process indexed by $V^2$ satisfies the following:
\begin{Theorem} \label{thm:Men-psi-1-Bernoulli}
There exist absolute constants $c_1$, $c_2$ and $c_3$ for which the following holds. If $H$ is a class of functions, then for every $r,u \geq c_1$, with $\mu^N$-probability at least $1-2\exp(-c_2u)$, $V=P_\sigma H$ satisfies that
$$
\sup_{v \in V} \left|\sum_{i=1}^N \eps_i v_i^2 \right| \lesssim ru^2\left(\gamma_2(H,\psi_2)+\sqrt{N}\sup_{h \in H} \|h\|_{\psi_1}\right) \cdot \gamma_2(H,\psi_2)
$$
with probability at least $1-2\exp(-c_3r^2)$ with respect to the Bernoulli random variables.
\end{Theorem}

Theorem \ref{thm:Men-psi-1-Bernoulli} is rather restricted because the $\psi_2$-based complexity parameter seems too strong in many situations, as does the assumption that $H$ is a bounded subset of $L_{\psi_1}$. Here, we will try to impose as few assumptions as possible on $H$.

Let $H$ be a class of functions on $(\Omega,\mu)$. For every $u>0$ we will define three events in the product space $\Omega^N$, which will be denoted by $\Omega_{1,u}$, $\Omega_{2,u}$ and $\Omega_{3,u}$. On the event $\Omega_{1,u} \cap \Omega_{2,u} \cap \Omega_{3,u}$, the random set $P_\sigma H=\{(h(X_i))_{i=1}^N: h \in H\}$ will be well behaved for the right choice of functionals $\theta_{s}$ and $\phi$. We will then study cases in which the event $\Omega_{1,u} \cap \Omega_{2,u} \cap \Omega_{3,u}$ has high probability.

\begin{Definition} \label{def:s-0}
 For $(\eta_s)_{s \geq 0}$ as above, set $s_0 \geq 0$ to be the first integer for which $\eta_s \geq \log(eN)$.

For every $s \in \{s : \log(eN) \leq \eta_s \leq N\}$, let $\ell_s$ be the largest integer in $\{1,...,N\}$ for which $\eta_s \geq \ell \log(eN/\ell)$, and if $ \eta_s \leq \log(eN)$, set $\ell_s =1$.
 \end{Definition}

The motivation for this definition is the following. If $E_k$ is the collection of subsets of $\{1,...,N\}$ of cardinality $k$, $s_0$ is the level above which one may find $k$ for which the cardinalities $|E_k|$ and $|H_s|$ are comparable. Indeed, when $s<s_0$, $|E_1|$ can be significantly larger than $|H_s| =2^{\eta_s}$, but when $s \geq s_0$, $\log |H_s|$ and $\log |E_{\ell_s}|$ are of the same order, and thus one may simultaneously control every function in $H_s$ and every subset in $E_{\ell_s}$ at no extra price. The main idea of the proofs in this section is to try and balance these two quantities as much as possible.

Observe that since $(\eta_s)_{s=0}^\infty$ grows exponentially, so does $(\ell_s)_{s \geq s_0}$.

\begin{Definition} \label{Omega-u-1}
 For an admissible sequence $(H_s)_{s \geq 0}$ and a sequence of functionals $(\theta_{u,s})_{s \geq s_0}$, let $\Omega_{1,u}$ be the event for which, for every $h \in H$, the following holds:
\begin{description}
\item{1.} for every $\log (eN) \leq \eta_s \leq N$, $\left(\sum_{i=1}^{u\ell_{s+1}} ((\Delta_s h)^2(X_i))^*\right)^{1/2} \leq \theta_{u,s}(\Delta_s h)$, (and if the $u\ell_{s+1} > N$ then the sum terminates at $N$).
\item{2.} for every $\eta_s > N$, $\left(\sum_{i=1}^{N} ((\Delta_s h)^2(X_i))^*\right)^{1/2} \leq \theta_{u,s}(\Delta_s h)$.
\item{3.} $\left(\sum_{i=1}^{u\ell_{s_0+1}} ((\pi_{s_0} h)^2(X_i))^* \right)^{1/2} \leq \theta_{u,s_0}(\pi_{s_0}h)$.
\end{description}
\end{Definition}
The set $\Omega_{1,u}$ is the subset of $\Omega^N$ in which the functionals $\theta_{u,s}$ yield a good bound on the $\ell_2$ norm of the ``relatively large" coordinates of each increment when $s \geq s_0$. In contrast, on the set $\Omega_{2,u}$ the smaller coordinates will be controlled for $s \geq s_0$. One of the key points of the proof is finding an estimate on the $\ell_2^n$ norm on these coordinates, but doing so without any real concentration phenomenon for sums of i.i.d. random variables coming to one's aid.

Formally, to define the set $\Omega_{2,u}$, first fix a random variable $Y$, an integer $N$ and $\eps>0$. For every $j \leq N$ let $\delta_j=(j/eN)^{(1+\eps)}$, set
$$
y_j=\inf\{y: Pr(|Y| \geq y_j) \leq \delta_j\},
$$
and without loss of generality, we will assume that the infimum is attained.

For every $1 \leq k \leq N$, let
$$
f_u(Y,k)=\kappa_3 \sqrt{u} \left(\sum_{\{j : 2^j \leq \lceil k/u \rceil\}} 2^j y^2_{2^j} \right)^{1/2},
$$
where $\kappa_3$ is a suitable chosen absolute constant.

The motivation for this definition is the following observation, showing that with high probability, the ``tail" of a sum of i.i.d random variables can be controlled using $f$.

\begin{Lemma} \label{lemma:tail-selectors}
There exist absolute constants $c_1$ and $c_2$ for which the following holds. For every integer $\ell$ and $u \geq c_1/\eps$, with probability at least $1-2\exp(-c_2u\eps \ell\log(eN/\ell))$, for every integer $k>u\ell$,
$$
\sum_{i=u\ell+1}^{k} (Y_i^2)^*  \leq f_u(Y,k)
$$
\end{Lemma}

\proof
Since $Pr(|Y| \geq y_j) \leq \delta_j =(j/eN)^{1+\eps}$ then for $u \geq 1$,
\begin{align*}
Pr(Y^*_{uj} \geq y_j) \leq & \binom{N}{uj} \delta_j^{uj} \leq \exp(uj\log(eN/uj)-(1+\eps)uj \log(eN/j))
\\
\leq & \exp(- \eps uj \log(eN/j)).
\end{align*}
Thus, summing over $\{j=\lceil \ell+2^i/u \rceil: \ 2^i \leq k-u\ell \}$, it follows that with probability at least $1-\exp(-c_1 \eps u \ell \log(eN/\ell))$, if $2^i \leq k-u\ell$ then $Y^*_{u\ell+2^i} \leq y_{\lceil \ell+2^i/u \rceil}$. Therefore,
\begin{align*}
\sum_{j=u\ell+1}^{k} (Y_j^2)^* \leq & \sum_{\{i : 2^i \leq k-u\ell\}} 2^i (Y^2_{u\ell+2^{i-1}})^* \leq \sum_{\{i : 2^i \leq k-u\ell\}} 2^i y_{\ell+2^{i-1}/u}^2
\\
\leq & c_2u \sum_{\{j : 2^j \leq k/u\}} 2^i y_{2^i}^2 = f^2(Y,k),
\end{align*}
where the last inequality is evident by a change of variables.
\endproof

We will also need the following ``global" counterpart of the functional $f$.
\begin{Definition} \label{def:global-parameters}
Given a class of functions $H$, an integer $N$ and $\eps>0$, set
$$
z_j=\inf\{z:\sup_{h \in H} Pr(|h| \geq z_j) \leq (j/eN)^{1+\eps}\}.
$$
For every $k \leq N$ and $u \geq 1$, let
$$
F_u(k)=\kappa_3 \sqrt{u}\left( \sum_{\{j : 2^j \leq k/u\}} 2^j z_{2^j}^2 \right)^{1/2}
$$
\end{Definition}
Clearly, for every $h \in H$ and every $k$, $f_u(h,k) \leq F_u(k)$.

\begin{Definition} \label{Def:global-F}
Let $\Omega_{2,u}$ be the event on which, for every $h \in H$, every $s \geq s_0$ and every $j > u \ell_s$
\begin{description}
\item{1.} $\left(\sum_{i=u\ell_s+1}^j ((\Delta_s h)^2(X_i))^*\right)^{1/2} \leq f_{u}(\Delta_s h, j)$,
\item{2.} $\left(\sum_{i=u\ell_s+1}^j ((\pi_{s} h)^2(X_i))^*\right)^{1/2} \leq F_u(j)$.
\end{description}
\end{Definition}
The final set, $\Omega_{3,u}$ is very close in nature to $\Omega_{2,u}$. It is needed to control the coordinates of ``very small" increments -- when $s<s_0$, if such an integer exists.
\begin{Definition} \label{def:omega-3-u}
If $\eta_0<\log(eN)$, let $\Omega_{3,u}$ be the event on which for every $h \in H$, every $0 \leq s < s_0$ and $1 \leq j \leq N$,
$$
\left(\sum_{i=1}^j ((\Delta_s h)^2(X_i))^* \right)^{1/2} \leq f_u(\Delta_s h, j), \ \ \ \left(\sum_{i=1}^j ((\pi_{s_0} h)^2(X_i))^* \right)^{1/2} \leq F_u(j).
$$
If $\eta_0 \geq \log(eN)$ set $\Omega_{3,u}=\Omega^N$.
\end{Definition}

It turns out that on the event $\Omega_{1,u} \cap \Omega_{2,u} \cap \Omega_{3,u}$, the set $P_\sigma H$ is indeed well behaved. Let
\begin{equation} \label{eq:gamma-u}
\gamma_u= \inf \sup_{h \in H} \sum_{s >s_0} \theta_{u,s}(\Delta_s h),
\end{equation}
with the infimum is taken with respect to all $(\eta_s)$-admissible sequences. From here on we will assume that $(H_s)_{s \geq 0}$ is an almost optimal $(\eta_s)_{s \geq 0}$-admissible sequence.

\begin{Lemma} \label{lemma:on-the-good-event}
There exists absolute constants $c_1$ and $c_2$ for which the following holds. Let $(\theta_{u,s})_{s \geq s_0}$ be functionals, and for $s<s_0$ set $\theta_{u,s}=0$. For every $u \geq c_1$, on the event $\Omega_{1,u} \cap \Omega_{2,u} \cap \Omega_{3,u}$, for every $h \in H$ and $I \subset \{1,...,N\}$,
\begin{description}
\item{1.} if $\eta_s \leq N$ then
$$
\left(\sum_{i \in I} (\Delta_s h)^2(X_i)\right)^{1/2} \leq \theta_{u,s}(\Delta_s h)+ f_{u}(\Delta_s h,|I|),
$$
and if $\eta_s >N$,
$$
\left(\sum_{i \in I} (\Delta_s h)^2(X_i) \right)^{1/2} \leq \theta_{u,s}(\Delta_s h).
$$
\item{2.}
$$
\left(\sum_{i \in I} h^2(X_i) \right)^{1/2} \leq \gamma_u+\sum_{\{i: 2^i \leq |I|\}} F(c_2u2^i)+R_{s_0}(h,I),
$$
\end{description}
where $R_{s_0,I}(h)=\theta_{u,0}(\pi_{0}h)$ if $s_0=0$ and $R_{s_0}(h,I) =\min\{\theta_{u,s_0}(\pi_{s_0}h), F_u(|I|)\}$ otherwise.
\end{Lemma}

\proof
First, assume that $\log(eN) \leq \eta_s \leq N$ (i.e. $s \geq s_0$) and recall that $\ell_s$ is the largest integer for which $\eta_s \geq \ell \log(eN/\ell)$. If $|I| \leq u\ell_s$ then the claim follows from the definition of $\theta_{u,s}$ and the set $\Omega_{1,u}$. If $|I| \geq u \ell_s$, then
$$
\left(\sum_{i \in I} (\Delta_s h)^2(X_i)
 \right)^{1/2} \leq \left(\sum_{i=1}^{u \ell_s} ((\Delta_s h)^2(X_i))^*\right)^{1/2} + \left(\sum_{i=u\ell_s +1}^{|I|} ((\Delta_s h)^2(X_i))^*\right)^{1/2},
$$
and the claim is evident from the definition of the function $f_u$ and the set $\Omega_{2,u}$.

If, on the other hand, $\eta_s < \log(eN)$ then $s_0>0$ and the assertion follows from the definition of $\Omega_{3,u}$.

The second part of (1) follows from the definition of $\Omega_{1,u}$.

Turning to (2), we shall treat two cases. First, consider the case $|I| \geq u\ell_{s_0}$ and observe that it suffices to estimate $(\sum_{i=1}^{u\ell_{s+1}} ((\pi_s h)^2(X_i))^*)^{1/2}$. Indeed, let $s$ be an integer for which $u\ell_s \leq |I| < u\ell_{s+1}$. Since $\ell_{s+1}$ is nondecreasing, then on $\Omega_{1,u}$,
\begin{align*}
\left(\sum_{i=1}^{u\ell_{s+1}} (h^2(X_i))^*\right)^{1/2} \leq & \sum_{j \geq s+1} \left(\sum_{i=1}^{u\ell_{s+1}} ((\Delta_j h)^2(X_i))^*\right)^{1/2}+\left(\sum_{i=1}^{u\ell_{s+1}} ((\pi_s h)^2(X_i))^* \right)^{1/2}
\\
\leq & \sum_{j \geq s+1} \theta_{u,j}(\Delta_j h)+\left(\sum_{i=1}^{u\ell_{s+1}} ((\pi_s h)^2(X_i))^* \right)^{1/2}.
\end{align*}

If $J \subset I$ is the set of the largest $u\ell_s$ coordinates of $\left((\pi_{s} h)(X_i)\right)_{i=1}^N$ in $I$, then the coordinate projections satisfy that
$$
P_I ((\pi_{s} h)(X_i))_{i=1}^N= P_{J} ((\pi_{s-1} h)(X_i))_{i=1}^N +  P_{J} ((\Delta_{s} h)(X_i))_{i=1}^N + P_{I \backslash J} ((\pi_{s} h)(X_i))_{i=1}^N,
$$
and thus,
\begin{align*}
& \left(\sum_{i \in I} (\pi_{s} h)^2(X_i) \right)^{1/2} \leq \left(\sum_{i =1}^{u \ell_{s+1}} \left((\pi_{s} h)^2(X_i)\right)^* \right)^{1/2}
\\
\leq & \max_{|I_1|=u\ell_s} \left(\sum_{i \in I_1} (\pi_{s-1} h)^2(X_i) \right)^{1/2}
+ \max_{|I_1|=u\ell_s} \left(\sum_{i \in I_1} (\Delta_{s} h)^2(X_i) \right)^{1/2}
\\
+ & \left(\sum_{i=u\ell_s+1}^{u\ell_{s+1}} \left((\pi_s h)^2(X_i)\right)^* \right)^{1/2}.
\end{align*}
Hence, if we set ${\cal U}_{j,s}(h)=  \max_{|I|=u\ell_{j+1}} \left(\sum_{i \in I} (\pi_s h)^2(X_i) \right)^{1/2}$ then for every $s$,
and every $h \in H$
\begin{align*}
{\cal U}_{s,s}(h) \leq & {\cal U}_{s-1,s-1}(h) +  \max_{|I_1|=u\ell_{s}} \left(\sum_{i \in I_1} (\Delta_{s} h)^2(X_i) \right)^{1/2}
\\
+ & \left(\sum_{i=u\ell_s+1}^{u\ell_{s+1}} \left((\pi_{s} h)^2(X_i)\right)^* \right)^{1/2}
\\
\leq & {\cal U}_{s-1,s-1}(h) + \theta_{u,s}(\Delta_s h) + F_{u}(u\ell_{s+1}).
\end{align*}
Summing over all $s > s_0$,
$$
{\cal U}_{s,s}(h) \leq \sum_{j=s_0+1}^{s} \theta_{u,j}(\Delta_j h)+ \sum_{j=s_0+1}^{s+1} F_{u}(u\ell_{j})  +  {\cal U}_{s_0,s_0}(h),
$$
and thus, for every $h \in H$ and every $I \subset \{1,...,N\}$,
\begin{equation*}
\left(\sum_{i \in I} h^2(X_i)  \right)^{1/2}
\leq \sum_{s > s_0} \theta_{u,s}(\Delta_j h)+ \sum_{\{s > s_0 : \ell_s \leq |I| \}} F_{u}(u\ell_{s+1}) +  {\cal U}_{s_0,s_0}(h).
\end{equation*}
Next, one has to bound
$\sup_{h \in H} \max_{|I| \leq u\ell_{s_0+1}} \left(\sum_{i \in I} (\pi_{s_0} h)^2(X_i) \right)^{1/2}$. This is at most $\theta_{u,s_0}(\pi_{s_0}h)$ on $\Omega_{1,u}$ and when $s_0>0$, it is also bounded by $F_u(u\ell_{s_0}) \leq F_u(|I|)$ on $\Omega_{3,u}$.

The claim in this case follows since $\ell_s$ grows exponentially for $s \geq s_0$, and thus
$$
\sum_{\{s \geq s_0 :\ell_s \leq |I|\}} F_{u}(u\ell_{s+1}) \leq \sum_{\{i : 2^i \leq |I|\}} F_{u}(cu2^i)
$$
for a suitable absolute constant $c$.

Turning to the second case, if $|I| \leq u\ell_{s_0}$, note that
\begin{align*}
\left(\sum_{i \in I} h^2(X_i)\right)^{1/2} \leq & \sum_{s>s_0} \left(\sum_{i \in I} (\Delta_s h)^2(X_i)\right)^{1/2} + \left(\sum_{i \in I} (\pi_{s_0} h)^2(X_i)\right)^{1/2}
\\
\leq & \sum_{s>s_0} \theta_{u,s}(\Delta_s h) + \min\{\theta_{u,s_0}(\pi_{s_0}h),F_u(|I|)\}.
\end{align*}
\endproof

For Lemma \ref{lemma:on-the-good-event} to have any meaning, one has to identify the functionals $f_u$, $F_u$ and $\theta_{u,s}$ in the cases one is interested in. Our next goal is to study the functions $f_u$ and $F_u$ under various tail assumptions on functions in $H$, and naturally, the two families of tail estimates we will be interested in are when $H$ has a bounded diameter in $L_{\psi_\beta}$ or in $L_q$ for $q>2$.

If $H \subset L_{\psi_\beta}$, then for every $h \in H$, $Pr(|h| \geq y) \leq \exp(-(y/\|h\|_{\psi_\beta})^\beta)$. Thus, for $\eps \geq 1$ and every $j$,
$$
y_j \lesssim \eps \|h\|_{\psi_\beta} \log^{1/\beta}(eN/j), \ \ \ z_j \lesssim \eps \sup_{h \in H} \|h\|_{\psi_\beta} \log^{1/\beta}(eN/j).
$$
Hence, if $d_{\psi_\beta}=\sup_{h \in H} \|h\|_{\psi_\beta}$, then
\begin{align*}
F_u(i) & \lesssim \eps \sqrt{u}\left(\sum_{j=1}^{\log_2 i} 2^j  z_{2^j}^2 \right)^{1/2} \lesssim_\beta \eps \sqrt{u} d_{\psi_\beta} \left(\sum_{j=1}^{\log_2 i} 2^j \log^{2/\beta}(eN/2^j) \right)^{1/2}
\\
& \lesssim_\beta \eps \sqrt{u}d_{\psi_\beta} \sqrt{i} \log^{1/\beta}(eN/i) \sim_\beta \eps \sqrt{u}d_{\psi_\beta}\phi_\beta(i),
\end{align*}
and in a similar fashion,
\begin{equation*}
f_{u}(h,i) \lesssim_\beta \eps \sqrt{u}\|h\|_{\psi_\beta} \sqrt{i} \log^{1/\beta}(eN/i) \sim_\beta \eps \sqrt{u}\|h\|_{\psi_\beta}\phi_\beta(i).
\end{equation*}

Using the same argument, if $h \in L_q$ then
$Pr(|h| \geq \|h\|_{L_q}y) \leq 1/y^q$ and for any $0 < \eps < q/2 -1$, $y_j =\|h\|_{L_q}(N/j)^{(1+\eps)/q}$. If
$\sup_{h \in H} \|h\|_{L_q} = d_{L_q}$, $q>2$ and $c_{q,\eps}=1-2(1+\eps)/q$ then
\begin{align*}
F_u(i) \lesssim & \sqrt{u}d_{L_q} \left(\sum_{j=1}^{\log_2 i} 2^j (N/2^j)^{2(1+\eps)/q} \right)^{1/2} \lesssim c^{-1}_{q,\eps}\sqrt{u}d_{L_q} \sqrt{i}\left(\frac{N}{i}\right)^{(1+\eps)/q}
\\
\sim & c_{q,\eps}^{-1}\sqrt{u}d_{L_q}\phi_{q,\eps}(i),
\end{align*}
and
$$
f_u(h,i) \lesssim c^{-1}_{q,\eps} \sqrt{u}\|h\|_{L_q}  \sqrt{i}\left(\frac{N}{i}\right)^{(1+\eps)/q} \sim c_{q,\eps}^{-1} \sqrt{u}\|h\|_{L_q}\phi_{q,\eps}(i).
$$
Combining these observations with the estimates of Lemma \ref{lemma:on-the-good-event} and noting that if $s_0>0$ then $R_{s_0,I}(h) \lesssim \sqrt{u}d_{L_q} \phi_{q,\eps}(|I|)$, one reaches the following corollary.
\begin{Corollary} \label{cor:middle-step}
Let $(\theta_{u,s})_{s \geq s_0}$ be a sequence of functionals and for $s < s_0$ let $\theta_{u,s}=0$. If $H$ is bounded in $L_q$ for $q>2$, then
on $\Omega_{1,u} \cap \Omega_{2,u} \cap \Omega_{3,u}$, for every $h \in H$ and every $I \subset \{1,...,N\}$
\begin{description}
\item{1.} if $\eta_s \leq N$,
$$
\left(\sum_{i \in I} (\Delta_s h)^2(X_i) \right)^{1/2} \lesssim \theta_{u,s}(\Delta_s h)+ c_{q,\eps}^{-1} \sqrt{u}\|\Delta_s h\|_{L_q}\phi_{q,\eps}(|I|),
$$
and if $\eta_s > N$ then
$$
\left(\sum_{i \in I} (\Delta_s h)^2(X_i) \right)^{1/2} \lesssim \theta_{u,s}(\Delta_s h).
$$
\item{2.}
$$
\left(\sum_{i \in I} h^2(X_i) \right)^{1/2} \lesssim \sum_{s \geq 0} \theta_{u,s}(\Delta_s h)+ c_{q,\eps}^{-1} \sqrt{u}d_{L_q}\phi_{q,\eps}(|I|).
$$
\end{description}
A similar bound holds when $H$ is bounded in $L_{\psi_\beta}$.
\end{Corollary}

\section{Estimates on $\Omega_{i,u}$ and the choice of functionals} \label{sec:estimates}
We will begin by showing that $\Omega_{2,u}$ is a large set, almost regardless of any assumptions on $\phi$, an observation that is based on the same idea as Lemma \ref{lemma:tail-selectors}.
\begin{Lemma} \label{lemma:est-Omega-2,u}
There exist absolute constants $c_1$ and $c_2$ such that, for every $\eps>0$ and $u \geq c_1/\eps$,  $Pr(\Omega_{2,u}) \geq 1-2\exp(-c_2\eps u \eta_{s_0})$.
\end{Lemma}

\proof
Recall that by Lemma \ref{lemma:tail-selectors}, for any random variable $Y$, with probability at least $1-2\exp(-c_1u\eps \ell\log(eN/\ell))$, for every integer $k>u\ell$,
\begin{equation} \label{eq:monotone-single}
\sum_{i=u\ell+1}^{k} (Y_i^2)^*  \leq f_u(Y,k).
\end{equation}
 Let $\ell=\ell_{s_0}$, and since $\eta_s \sim \ell_s \log(eN/\ell_s)$ and $|\Delta_s H| \lesssim 2^{\eta_s}$, then for $u \geq c_3/\eps$, \eqref{eq:monotone-single} holds uniformly for every $h \in \Delta_s H$
 with probability at least $1-2\exp(-c_4u \eps \eta_s)$. The analogous claim holds for functions in $H_s$ as well, with the uniform bound of $F_u$ replacing $f_u$. Summing over all $s \geq s_0$ and since $(\eta_s)$ grows exponentially, the claim follows.
\endproof

Since $\Omega_{2,u}$ is always large, and since $\Omega_{3,u}$ will behave in a very similar way when $s_0>0$, the crucial point in the construction of a good decomposition of $P_\sigma H$ is a correct choice of $\theta_{u,s}$ and estimates on $\Omega_{1,u}$.

The functionals $\theta_{u,s}$ capture the geometry of $H$, and thus have to be selected according to the information one has on the class. We will present two examples of such choices, each leading to one of our two main results. The first one will be based on ``global" structure like metric entropy, while the second uses accurate estimates on each ``chain".

\subsection{The ball  $B_2^n$ -- global estimates}
Let $\mu$ be an unconditional measure on $\R^n$, set $H=\{ \inr{t,\cdot} : t \in B_2^n\}$ to be a class of linear functionals on $(\R^n,\mu)$ -- and from here on we will identify the class $\{\inr{t,\cdot} : t \in T\}$ with its indexing set $T$. We will also assume that $\mu$ satisfies the $p$-small diameter, $L_q$ moment assumption for some $p>2$ and $q>2$; that is, $\mu$ is supported in $\kappa_1 n^{1/p} B_p^n$, and for every $x \in \R^n$, $\|\inr{x,\cdot}\|_{L_q} \leq \kappa_2 \|x\|_{\ell_2^n}$.

Let $\kappa_4 \geq 10$ be an absolute constant to be fixed later, set $2^{s_1} \sim n^{\delta}$ for $\delta < 1/2-1/2(p-1)$, and put
$$
\eta_s = \kappa_4 2^{s+s_1}\max\{\log(en/2^{s+s_1}),1\}.
$$
Note that $s_0 =0$ as long as $\eta_0 \sim 2^{s_1}\log(en/2^{s_1}) \geq \log(eN)$, i.e., if $n^{\delta} \log(n) \gtrsim \log (eN)$ - which we will assume is the case, since our main interest in when $N \sim n$.

If $X=(x_1,...,x_n)$ is distributed according to $\mu$ then for every $1 \leq \ell \leq n$,
set $M_\ell=\|(\sum_{i=1}^\ell (x_i^2)^* )^{1/2}\|_{L_\infty}$. Define the following functionals (which, in this case, will be constants depending only on $u$ and $s$): let $
\theta_{u,0} = c \sqrt{u} \eta_0^{1/2} n^{1/p} 2^{(s+s_1)(1/2-1/p)},
$
if $2^{s+s_1} \leq n$, set
$
\theta_{u,s} = c\sqrt{u}\eta_s^{1/2} n^{1/p} 2^{-(s+s_1)/p}
$
and if $\eta_s \geq n$ put $\theta_{u,s} = c \sqrt{u}\eta_s^{1/2} 2^{-2^s/n}$, where $c=c(\kappa_1,p,\delta)$.
\begin{Theorem} \label{thm:bounded-uncon}
For every $\kappa_1$, $p>2$ and $\delta<1/2-1/2(p-1)$ there exist constants $c_1, c_2$ and $c_3$ that depend only on $\kappa_1$, $p$ and $\delta$ for which the following holds. There is an $(\eta_s)_{s \geq 0}$-admissible sequence of $B_2^n$, for which, if $u \geq c_1$, then $Pr(\Omega_{1,u}) \geq 1-\exp(-c_2n^{\delta})$ and
$$
\sup_{t \in B_2^n} \sum_{s \geq 0} \theta_{u,s}(\inr{\Delta_s t,\cdot}) \leq c_3 \sqrt{u}\sqrt{n}.
$$
\end{Theorem}
Observe that by the $p$-small diameter assumption, $M_\ell \lesssim_p n^{1/p} \ell^{1/2 -1/p}$. Also, since $\mu$ is unconditional, then for every $I \subset \{1,...,n\}$ and $v$ supported on $I$,
\begin{equation} \label{eq:psi-2-l-infty}
\|\inr{v,\cdot}\|_{\psi_2} \lesssim \|v\|_{\ell_\infty^I} M_{|I|}.
\end{equation}
Indeed, by the unconditionality of $\mu$, $(x_1,...,x_n)$ has the same distribution as $(\eps_1x_1,...,\eps_nx_n)$. Hence, for every $r \geq 1$
\begin{align*}
\|\inr{v,\cdot}\|_{L_r} \sim & (\E_X \E_\eps |\sum_{i \in I} \eps_i x_i v_i|^r)^{1/r} \lesssim \left(\E_X r^{r/2}(\sum_{i \in I} v_i^2 x_i^2)^{r/2}\right)^{1/r}
\\
\lesssim & \sqrt{r}\|v\|_{\ell_\infty^I} M_{|I|}.
\end{align*}
We will also need a few $\psi_2$ entropy estimates. Set $B_{\psi_2}=\{v \in \R^n : \|\inr{v,\cdot}\|_{\psi_2} \leq 1\}$, and for $K,L \subset \R^n$ denote by $N(K,L)$ the minimal number of translates of $L$ needed to cover $K$.
\begin{Lemma} \label{lemma:entropy-subset}
If $I \subset \{1,...,n\}$ then for every $\eps>0$, $\log N(B_2^I,\eps B_{\psi_2}) \lesssim M_{|I|}^2/\eps^2$. Moreover, for $\eps \leq 1$, $\log N(B_2^n,\eps B_{\psi_2}) \lesssim n\log(2/\eps)$.
\end{Lemma}
\proof
By the dual Sudakov inequality (see, e.g. \cite{LT}), if $B_{\| \ \|}$ is a unit ball of a norm on $\R^I$ and $G=(g_i)_{i \in I}$ is a standard Gaussian vector on $\R^I$, then $\log N(B_2^n,\eps B_{\| \ \|}) \lesssim (\E \|G\|)^2/\eps^2$. Since $\|f\|_{\psi_2} \leq \E \exp(f^2)$ and $(\sum_{i \in I} x_i^2 )^{1/2} \leq M_{|I|}$ almost surely, then by changing the order of integration,
$$
\E \|G/cM_{|I|}\|_{\psi_2} \leq \E_X \E_G (\exp( (\sum_{i \in I} g_i x_i)^2/c^2M_{|I|}^2)|X) \leq 2
$$
for a suitable absolute constant $c$, proving the first part.

For the second part, note that $N(B_2^n,\eps B_{\psi_2}) \leq N(B_2^n, B_{\psi_2}) \cdot N(B_{\psi_2},\eps B_{\psi_2})$. By the first part, $\log N(B_2^n, B_{\psi_2}) \lesssim n$, while a standard volumetric estimate shows that $N(B_{\psi_2},\eps B_{\psi_2}) \leq (5/\eps)^n$.
\endproof
Next, let us define the sets $T_s$. If $2^{s+s_1} > n$, let $T_s$ be a maximal $\eps_s$ separated subset of $B_2^n$ relative to the $\psi_2$ norm and of cardinality $2^{\eta_s}$. If $2^{s+s_1} \leq n$, let $T_s$ be a maximal $\eps_s$ separated subset of $U_{2^{s+s_1}}=\{x \in B_2^n : |{\rm supp}(x)| \leq 2^{s+s_1}\}$ with respect to the $\psi_2$ norm, and of cardinality $2^{\eta_s}$. Given a vector $t \in B_2^n$, we will define the functions $\pi_s$ as follows. If $2^{s+s_1} > n$, $\pi_s t$ is a best $\psi_2$ approximation of $t$ in $T_s$.
For $2^{s+s_1} \leq n$ one combines approximation and dimension reduction. Set $s_*$ to satisfy that $2^{s_*+s_1}=n$ (and without loss of generality we will assume that such an integer exists). If $v=\pi_{s_*} t$, let $I_{n/2}$ be the set of the largest $n/2$ coordinates of $v$, and put $\pi_{s_*-1} t$ to be the best approximation of the coordinate projection $P_{I_{n/2}}v$ in $T_{s_*-1}$, and so on.
\begin{Lemma} \label{lemma:psi-2-estimate}
There exists an absolute constant $c$ such that for every $t \in B_2^n$, if $s > s_*$ (i.e., if $\eta_s > \kappa_3 n$), then
$\|\inr{\Delta_s t,\cdot}\|_{\psi_2} \leq c2^{-2^{s+s_1}/n}$, and if $0<s \leq s_*$ then $\|\inr{\Delta_s t,\cdot}\|_{\psi_2} \leq c2^{-(s+s_1)/2}M_{2^{s+s_1}}$.
\end{Lemma}
\proof
First consider $s > s_*$. Note that $\|\inr{\Delta_s t,\cdot}\|_{\psi_2} \leq \|\inr{t-\pi_{s} t,\cdot}\|_{\psi_2} + \|\inr{t-\pi_{s-1} t,\cdot}\|_{\psi_2} \leq \eps_{s} + \eps_{s-1}$, and by the covering numbers estimate from Lemma \ref{lemma:entropy-subset}, in that range $\eps_s \lesssim 2^{-2^{s+s_1}/n}$.

In the range $s \leq s_*$, $\Delta_s t = u+w$, where $w$ consists of the smallest $2^{s+s_1-1}$ coordinates of $\pi_{s} t \in B_2^I$ for some $|I|=2^{s+s_1}$, and $u$ is an $\eps_{s-1}$-approximation of the largest $2^{s+s_1-1}$ coordinates of $\pi_{s}t$. Therefore, $\|\inr{\Delta_s t, \cdot}\|_{\psi_2} \leq \|\inr{w,\cdot}\|_{\psi_2}+\eps_{s-1}$. Recall that for every such $s$,  $U_{2^{s+s_1}}$ is a union of $\binom{n}{2^{s+s_1}}$ balls of dimension $2^{s+s_1}$, then
\begin{align*}
\log N(U_{2^{s+s_1}},\eps B_{\psi_2}) \leq & 2^{s+s_1}\log(en/2^{s+s_1})+
\max_{|J|=2^{s+s_1}} \log N(B_2^J,\eps B_{\psi_2})
\\
\lesssim & 2^{s+s_1}\log(en/2^{s+s_1})+ M_{2^{s+s_1}}^2/\eps^2.
\end{align*}
Note that for a suitable choice of $\kappa_4$, $\log |T_s| \geq 2 \cdot 2^{s+s_1}\log(en/2^{s+s_1})$. Therefore, $\eps_s \leq 2^{-(s+s_1)/2}M_{2^{s+s_1}}$, and applying \eqref{eq:psi-2-l-infty}, $\|\inr{w,\cdot}\|_{\psi_2} \lesssim \|w\|_{\ell_\infty^I} M_{|I|} \lesssim 2^{-(s+s_1)/2}M_{2^{s+s_1}}$.
\endproof

\noindent {\bf Proof of Theorem \ref{thm:bounded-uncon}.} Observe that $\|Y\|_{\psi_2}^2=\|Y^2\|_{\psi_1}$, and thus, by a standard application of Bernstein's inequality, for every integer $m$,
$$
Pr\left( \sum_{i=1}^m Y_i^2 \geq m\|Y\|_{\psi_2}^2 t^2 \right) \leq 2\exp(-cm \min\{t^2,t^4\}).
$$
Therefore, if $w$ is large enough, then
\begin{align*}
& Pr \left( \sum_{i=1}^{u\ell_s} (Y_i^2)^* \geq w^2 \|Y\|_{\psi_2}^2 \cdot u\ell_s \log(eN/u\ell_s) \right)
\\
\leq & \binom{N}{u\ell_s} \cdot 2\exp(-cw^2 u\ell_s \log(eN/u\ell_s))
\leq 2\exp(-c_1w^2u\ell_s \log(eN/u\ell_s)).
\end{align*}
Moreover, $u\ell_s \log(eN/u\ell_s) \lesssim u\eta_s$, and for $u \geq 1$, $u\ell_s \log(eN/u\ell_s) \gtrsim \eta_s$, implying that with probability at least $1-2\exp(-c_2w^2\eta_s)$,
$$
\left(\sum_{i=1}^{u\ell_s} (Y_i^2)^*\right)^{1/2} \lesssim w u^{1/2} \|Y\|_{\psi_2} \eta_s^{1/2}.
$$
Also, with probability at least $1-2\exp(-c_2w^2\eta_s)$, if $\eta_s \geq N$ then
$$
\left(\sum_{i=1}^{N} (Y_i^2)^*\right)^{1/2} \lesssim w u^{1/2} \|Y\|_{\psi_2} \eta_s^{1/2}.
$$
Using Lemma \ref{lemma:psi-2-estimate} and summing the probability estimates, it is evident that with probability at least $1-2\exp(-c_3w^2\eta_0)$, the following holds: if $\eta_s \geq N$ then
$$
\sup_{t \in B_2^n} \left(\sum_{i=1}^N (\inr{\Delta_s t, X_i}^2)^*\right)^{1/2} \lesssim wu^{1/2} \eta_s^{1/2} 2^{-2^{s+s_1}/n},
$$
if $\kappa_4 n \leq \eta_s < N$, then
$$
\sup_{t \in B_2^n} \left(\sum_{i=1}^{u\ell_s} (\inr{\Delta_s t, X_i}^2)^*\right)^{1/2} \lesssim wu^{1/2} \eta_s^{1/2} 2^{-2^{s+s_1}/n},
$$
and if $s>0$ and $\eta_s \leq \kappa_4 n$ then
$$
\sup_{t \in B_2^n} \left(\sum_{i=1}^{u\ell_s} (\inr{\Delta_s t, X_i}^2)^*\right)^{1/2} \lesssim wu^{1/2} \eta_s^{1/2} 2^{-(s+s_1)/2}M_{2^{s+s_1}}.
$$
Finally, since $\eta_0=\kappa_4 2^{s_1}\log(en/2^{s_1})$ then $\ell_0 \lesssim 2^{s_1}$. Moreover, $|{\rm supp}(\pi_0 t)| \leq 2^{s_1}$ and by \eqref{eq:psi-2-l-infty}, $\|\inr{\pi_0 t ,\cdot}\|_{\psi_2} \leq  M_{2^{s_1}}$. Hence, with probability at least $1-2\exp(-c_4 w^2 \eta_0)$,
$$
\sup_{t \in B_2^n} \left(\sum_{i=1}^{u\ell_s} (\inr{\pi_0 t, X_i}^2)^*\right)^{1/2} \lesssim wu^{1/2} \eta_0^{1/2} M_{2^{s_1}}.
$$
Since $M_\ell \lesssim_p n^{1/p}\ell^{1/2-1/p}$, then $Pr(\Omega_{1,u}) \geq 1-2\exp(-c_5\eta_0)=1-2\exp(-c_52^{s_1})$ for the desired functionals $\theta_{u,s}$. It remains to choose $s_1$ and estimate $\sum_{s \geq 0} \theta_{u,s}$.

Note that if $2^{s_1} \sim n^{\delta}$ for $\delta < 1/2-1/2(p-1)$, then
\begin{equation} \label{eq:ball-theta-0}
\theta_0 \sim \eta_0^{1/2} M_{2^{s_1}} \sim_{\kappa_1,p} 2^{s_1/2}\log^{1/2}(en/2^{s_1}) n^{1/p} 2^{s_1(1/2-1/p)} \leq c_6(\kappa_1,p,\delta) \sqrt{n}.
\end{equation}
Also,
\begin{align}  \label{eq:ball-theta-small}
& \sum_{\{s>0:\eta_s \leq \kappa_3 n\}} \theta_s \lesssim \sum_{\{s: \eta_s \leq \kappa_3 n\}} \eta_s^{1/2} 2^{-(s+s_1)/2}M_{2^{s+s_1}} \nonumber
\\
\lesssim_{\kappa_1,p,\delta} & n^{1/p} \sum_{\{s: 2^{s+s_1} \leq n\}} 2^{(s_1+s)(1/2-1/p)}\log^{1/2}(en/2^{s+s_1}) \leq c_6(\kappa_1,p,\delta) \sqrt{n},
\end{align}
and
\begin{equation}  \label{eq:ball-theta-large}
\sum_{\{s: \eta_s > \kappa_3 n\}} \theta_s \lesssim_{\kappa_1,\kappa_3,p} \sum_{\{s:\eta_s > \kappa_3 n\}} 2^{(s+s_1)/2} 2^{-2^{s+s_1}/n} \leq c_6(\kappa_1,p,\delta) \sqrt{n}.
\end{equation}
\endproof

\begin{Corollary} \label{cor:decomp-sphere}
There exist absolute constants $c_1$, $c_2$ and $c_3$ and $c_4$ that depend on $\kappa_1,\kappa_2,p,\delta$, for which the following holds.
If $\mu$ is as above and $\eps>0$, then $B_2^n$ has an $(\eta_s)_{s \geq 0}$-admissible sequence $(T_s)_{s \geq 0}$ for which, for $u \geq c_1/\eps$, with probability at least $1-2\exp(-c_2\eps u n)-2\exp(-c_3n^{\delta})$, for every $t \in B_2^n$ and every $I \subset \{1,...,N\}$,
\begin{description}
\item{1.} if $\eta_s \leq N$,
$$
\left(\sum_{i \in I} (\inr{\Delta_s t,X_i})^2 \right)^{1/2} \leq  c_4 \theta_{u,s}+ c_{q,\eps}^{-1} \sqrt{u}\|\Delta_s t\|_{\ell_2^n} \phi_{q,\eps}(|I|),
$$
and if $\eta_s > N$ then
$$
\left(\sum_{i \in I} (\inr{\Delta_s t,X_i})^2 \right)^{1/2} \leq c_4 \theta_{u,s}.
$$
\item{2.}
\begin{align*}
\left(\sum_{i \in I} (\inr{t,X_i})^2 \right)^{1/2} \leq & c_4 \sum_s \theta_{u,s}+ c_{q,\eps}^{-1} \sqrt{u}\phi_{q,\eps}(|I|)
\\
\lesssim & \sqrt{u}\sqrt{n} + c_{q,\eps}^{-1} \sqrt{u}\phi_{q,\eps}(|I|).
\end{align*}
\end{description}
\end{Corollary}
We will separate our treatment to the cases $q>4$ and $2<q \leq 4$. First, if $q>4$, let $\eps=(q/4-1)/2$ and note that $c_{q,\eps} \geq 1/2$. Also, since $\|t\|_{\ell_2^n} \sim_{\kappa_2} \|\inr{t,\cdot}\|_{L_q} \lesssim_{\kappa_2} \|\inr{t,\cdot}\|_{\psi_2}$, then by the same computation as in \eqref{eq:ball-theta-0}, \eqref{eq:ball-theta-small} and \eqref{eq:ball-theta-large},
$$
B_4 = \sup_{t \in B_2^n} \sum_{s \geq 0} \eta_s^{1/2}\|\Delta_s t\|_{L_q} \lesssim_{\kappa_1,\kappa_2,p,\delta} \sqrt{n}.
$$
We thus have:
\begin{Theorem} \label{thm:Bernoulli-sphere-large-q}
For every $\kappa_1$, $\kappa_2$, $q>4$, $p>2$ and $\delta < 1/2-1/2(p-1)$, there exist constants $c_0$, $c_1$, $c_2$ and $c_3$ which depend on $\kappa_1$, $\kappa_2$, $p$, $q$ and $\delta$, and an absolute constant $c_4$ for which the following holds. If $\mu$ is as above, and $N \leq \exp(c_0n^{\delta})$, then for every $u \geq c_1$,
with $\mu^N$-probability at least $1-2\exp(-c_2 n^{\delta})$, $P_\sigma(B_2^n)$ satisfies that
$$
\sup_{t \in B_2^n} \left|\frac{1}{N}\sum_{i=1}^N \eps_i \inr{X_i,t}^2 \right| \leq c_3 ru\left(\sqrt{\frac{n}{N}}+\frac{n}{N}\right),
$$
with probability at least $1-2\exp(-c_4nr^2)$ relative to the Bernoulli random variables.
\end{Theorem}
Turning to the case $2<q \leq 4$, recall that for $0<\eps<q/2-1$, $B_{q,\eps}=\sum_{\{s : \eta_s \leq N\}} \eta_s^{1-2(1+\eps)/q} \|\Delta_s v\|_{L_q}$. Assume that $\mu$ is as above and satisfies the $p$-small diameter assumption for $p>q/(q/2-1)$. Then, for $0<\eps<q/2-1-q/p$ (i.e. if $1-(2(1+\eps)/q)-1/p>0$),
\begin{align*}
B_{q,\eps} & \lesssim \sum_{\{s : 2^{s+s_1} \leq n\}} (2^{s+s_1}\log(en/2^{s+s_1}))^{1-2(1+\eps)/q} 2^{-(s+s_1)/2} n^{1/p}2^{(s+s_1)(1/2-1/p)}
\\
& + \sum_{\{s : 2^{s+s_1} > n\}} 2^{(s+s_1)(1-2(1+\eps)/q)} 2^{-2^{(s+s_1)/n}}
\lesssim \frac{n^{1-2(1+\eps)/q}}{1-2(1+\eps)/q-1/p}.
\end{align*}
Therefore, one has
\begin{Theorem} \label{thm-Bernoulli-sphere-2<q<4}
Let $2< q \leq 4$, $p>(1-2/q)^{-1}$ and $0<\eps<2/q-1-q/p$. If $\mu$ and $\delta$ are as above, $u \gtrsim 1/\eps$, and $N \lesssim \exp(c_0n^{\delta})$, then
with $\mu^N$ probability at least $1-2\exp(-c_1 \eps u n)-2\exp(-c_2n^{\delta})$, $P_\sigma(B_2^n)$ satisfies that
$$
\sup_{t \in B_2^n} \left|\frac{1}{N}\sum_{i=1}^N \eps_i \inr{X_i,t}^2 \right| \lesssim_{\kappa_1,\kappa_2,p,\delta} \frac{ru}{(1-2(1+\eps)/q)^2} \left(\left(\frac{n}{N}\right)^{1-(2/q)-2\eps/q}+\sqrt{\frac{n}{N}}+\frac{n}{N}\right).
$$
with probability at least $1-2\exp(-c_3nr^2)$ relative to the Bernoulli random variables.

In particular, taking $\eps \sim 1/\log(eN/n)$, then for every such $N$ satisfying that $N \gtrsim_{\kappa_1,\kappa_2,q,p} n$, and any $u \geq_{\kappa_1,\kappa_2,q,p} \log(eN/n)$,
$$
\sup_{t \in B_2^n} \left|\frac{1}{N}\sum_{i=1}^N \eps_i \inr{X_i,t}^2 \right| \lesssim_{\kappa_1,\kappa_2,q,p} ru \left(\frac{n}{N}\right)^{1-2/q}.
$$
\end{Theorem}

\subsection{Unconditional log-concave measures} \label{sec:log-concave}
We will now present a different way of bounding $\Omega_{1,u}$ (and $\Omega_{3,u}$ if needed) by estimating the moments of the increments $\Delta_s h$, and selecting the functionals $\theta_{u,s}$ accordingly.

For every $s \geq 0$ and $h \in H$, set
$$
Z_s^2(h)=\sum_{i=1}^{\min\{u\ell_{s+1},N\}} ((\Delta_s h)^2(X_i))^*, \ \ \ Z_{s_0}^2(h)=\sum_{i=1}^{\min\{u\ell_{s_0+1},N\}} ((\pi_{s_0} h)^2(X_i))^*.
$$
In light of Theorem B, we will assume that $H$ is a bounded subset of $L_{\psi_1}$ (although what we do here can be extended to other moment assumptions), and thus one may control $\Omega_{2,u}$ using $\phi_{\beta}$ for $\beta=1$ and $\eps$ which will be selected later.

\begin{Lemma} \label{lemma:moment-idea}
There exist absolute constants $c_1$, $c_2$ and $c_3$ for which the following holds. For $u \geq c_1$, with probability at least $1-2\exp(-c_1u\eta_{s_0})$, for every $s \geq s_0$ and every $h \in H$,
$Z_s(h) \leq e\|Z_s(h)\|_{L_{2u\eta_{s+1}}}$.
\end{Lemma}
\proof If $Z$ is a nonnegative random variable then $Pr(Z \geq e\|Z\|_{L_q}) \leq \exp(-q)$. Thus, for a fixed $s$ and every $h \in H$,
$Z_s(h) \leq e\|Z_s(h)\|_{L_{2u\eta_{s+1}}}$ with probability at least $1-\exp(-2u\eta_{s+1})$. Since $\log |\Delta_s H| \lesssim \eta_{s+1}$ and because there are at most $\exp(u\ell_{s+1}\log(eN/u\ell_{s+1})) \leq \exp(u\eta_{s+1})$ subsets of $\{1,...,N\}$ of cardinality $u\ell_{s+1}$, the same probability estimate holds uniformly for every $h \in H$ (with a different constant). Summing the probabilities for every $s \geq s_0$ and repeating the same argument for $H_{s_0}$ concludes the proof.
\endproof

Next, one has to control the moments appearing in Lemma \ref{lemma:moment-idea}, which is based on the following result, due to Lata{\l}a \cite{Lat-moments}.
\begin{Theorem} \label{thm:latala}
Let $X_1,...,X_m$ be independent, distributed according to a nonnegative random variable $X$. Then for every $p \geq 1$,
$$
\|\sum_{i=1}^m X_i\|_{L_p} \sim \left\{ \frac{p}{r}\left(\frac{m}{p}\right)^{1/r}\|X\|_{L_r} \ : \ \max\{1,p/m\} \leq r \leq p \right\}.
$$
\end{Theorem}

\begin{Definition} \label{def:moment-norm}
If $X$ is a random variable, for every $p \geq 1$ set
$$
\|X\|_{(p)} = \sup_{1 \leq q \leq p} \frac{\|X\|_{L_q}}{\sqrt{q}}.
$$
\end{Definition}
The $(p)$-norms are a local version of the $\psi_2$ norm, and clearly $\|X\|_{(p)} \lesssim \|X\|_{\psi_2}$. Using those norms one may obtain a more compact expression for the required moments.
\begin{Lemma} \label{lemma:local-estimate}
There exist an absolute constant $c$ such that for every $h \in H$, every $s>s_0$ and every $u>0$,
$$
\|Z_s(h)\|_{L_{2u\eta_{s+1}}} \leq c\sqrt{u}\eta_{s+1}^{1/2} \|\Delta_s h\|_{(2u\eta_{s+1})}
$$
and
$$
\|Z_{s_0}(h)\|_{L_{2u\eta_{s_0+1}}} \leq c\sqrt{u}\eta_{s_0+1}^{1/2} \|\pi_{s_0} h\|_{(2u\eta_{s_0+1})}.
$$
\end{Lemma}

\proof Let $Y_i=h(X_i)$ and observe that for every $m$, $\|(\sum_{i=1}^m Y_i^2 )^{1/2}\|_{L_p} = \|\sum_{i=1}^m Y_i^2\|_{L_{p/2}}^{1/2}$. Since $m=u\ell_{s+1}$ and $p=2u\eta_{s+1}$ then $p/2 \geq m$. Also, for every $r \leq p$, $\|Y\|_{L_r} \leq \sqrt{r}\|Y\|_{(p)}$, and applying Theorem \ref{thm:latala},
$$
\|\sum_{i=1}^m Y_i^2\|_{L_{p/2}} \lesssim \|Y^2\|_{(p/2)} \sup_{2 \leq r \leq p} \frac{p}{\sqrt{r}}\left(\frac{m}{p}\right)^{1/r} \lesssim \|Y^2\|_{(p/2)} \frac{p}{\sqrt{\log(p/m)}}.
$$
Hence, for our choice of $p$ and $m$,
$$
\|\sum_{i=1}^m Y_i^2\|_{L_{p/2}}^{1/2} \lesssim \sqrt{u}\eta_{s+1}^{1/2}\|Y^2\|^{1/2}_{(u\eta_{s+1})}
=\sqrt{u}\eta_{s+1}^{1/2}\|Y\|_{(2u\eta_{s+1})}.
$$
\endproof

\begin{Corollary} \label{cor:chaining-bound-local}
There exist absolute constants $c_1$ and $c_2$ for which the following holds. If, for $s > s_0$,
$$
\theta_{u,s}(\Delta_s h) =c_1\sqrt{u}\eta_{s+1}^{1/2}\|\Delta_s h\|_{(2u\eta_{s+1})}
$$
and
$$
\theta_{u,s_0}(\pi_{s_0} h) =c_1\sqrt{u}\eta_{s_0+1}^{1/2}\|\pi_{s_0} h\|_{(2u\eta_{s_0+1})},
$$
then $Pr(\Omega_{1,u}) \geq 1-2\exp(-c_2u\eta_{s_0})$.
\end{Corollary}

Next, assume that $s_0>0$, and thus one has to bound $Pr(\Omega_{3,u})$.
\begin{Lemma} \label{lemma:Omega-3-u-local}
There exists absolute constants $c_1$, $c_2$ and $c_3$ such that, for every $u \geq c_1$, with probability at least $1-2\exp(-c_2u \log N)$, for every $0 \leq s <s_0$ and every $h \in H$,
$$
\left(\sum_{i=1}^j ((\Delta_s h)^2(X_i))^* \right)^{1/2} \leq c_3 u \|\Delta_s h\|_{\psi_1} \sqrt{j}\log(eN/j) \sim u \|\Delta_s h\|_{\psi_1} \phi_1(j),
$$
and a similar bound holds for $\pi_{s_0} h$.
\end{Lemma}

\proof
Recall that for a fixed $\eps>0$ and every $i$, $Pr(Y_i^* \geq y_i ) \leq \exp(-\eps i \log(eN/i))$. Let $\eps \sim u \geq 1$ and observe that if $Y \in L_{\psi_1}$, then $y_i \lesssim u\|Y\|_{\psi_1}\log(eN/i)$ and
\begin{equation} \label{eq:monotone-i-local}
Pr(\exists i \leq N : Y_i^* \geq y_i) \leq \exp(-c_1u \log N).
\end{equation}
Since the cardinality of the set $\cup_{s < s_0} \Delta_s H$ is at most $\sum_{s < s_0} 2^{\eta_{s+1}} \lesssim N^{c_2}$, \eqref{eq:monotone-i-local} holds uniformly with probability at least $1-\exp(-c_3u\log N)$ for $u \geq c_4$. Therefore, on that event, for every $0\leq s <s_0$ and every $j$,
$$
\left(\sum_{i=1}^j ((\Delta_s h)^2(X_i))^* \right)^{1/2} \leq c_5 u \|\Delta_s h\|_{\psi_1} \sqrt{j}\log(eN/j).
$$
An identical argument holds for $\left(\sum_{i=1}^j ((\pi_{s_0} h)^2(X_i))^* \right)^{1/2}$.
\endproof
Therefore, the event $\Omega_{1,u} \cup \Omega_{2,u} \cup \Omega_{3,u}$ has high probability, leading to the following decomposition result.
\begin{Corollary} \label{cor:decomposition-psi-1}
There exist absolute constants $c_1$ and $c_2$ for which the following holds. For every $u \geq c_1$, with probability at least $1-2\exp(-c_2u\log N)$, for every $h \in H$ and every $I \subset \{1,...,N\}$,
\begin{description}
\item{1.} if $\eta_s \leq N$,
$$
\left(\sum_{i \in I} ((\Delta_s h)^2(X_i))^* \right)^{1/2} \lesssim \sqrt{u}\eta_{s+1}^{1/2}\|\Delta_s h\|_{(2u\eta_{s+1})}+  u\|\Delta_s h\|_{\psi_1} \phi_{1}(|I|),
$$
and if $\eta_s > N$ then
$$
\left(\sum_{i \in I} ((\Delta_s h)^2(X_i))^* \right)^{1/2} \lesssim \sqrt{u}\eta_{s+1}^{1/2}\|\Delta_s h\|_{(2u\eta_{s+1})}.
$$
\item{2.}
$$
\left(\sum_{i \in I} (h^2(X_i))^* \right)^{1/2} \lesssim \sqrt{u}\sum_{s>s_0}\eta_{s+1}^{1/2}\|\Delta_s h\|_{(2u\eta_{s+1})} + ud_{\psi_1} \phi_{1}(|I|)+R_{s_0}(h),
$$
\end{description}
where $R_{s_0}(h) \lesssim \sqrt{u}\eta_1\|\pi_0 h\|_{(2u\eta_1)}$ if $s_0=0$ and $R_{s_0}(h) \lesssim u d_{\psi_1}\phi_1(|I|)$ otherwise.
\end{Corollary}
\begin{Remark}
Note that $\|\Delta_s h\|_{(2u\eta_{s+1})} \lesssim \|\Delta_s h\|_{\psi_2}$, and thus one may take $\theta_{u,s} \sim \sqrt{u} \eta^{1/2}_{s+1}\|\Delta_s h\|_{\psi_2}$. If $\eta_0=0$ and $\eta_s =2^s$ for $s \geq 1$, then for an almost optimal admissible sequence,
$$
\sum_{s>s_0}\eta_{s+1}^{1/2}\|\Delta_s h\|_{(2u\eta_{s+1})} \lesssim \gamma_2(H,\psi_2).
$$
Although this estimate leads to an alternative proof of Theorem \ref{thm:Men-psi-1-Bernoulli}, it is not sharp enough to prove Theorem B, as the latter requires more accurate bounds on $\|\Delta_s h\|_{(2u\eta_{s+1})}$.
\end{Remark}
 From here on we will assume that $\eta_0=0$ and that $\eta_s =2^s$ for $s \geq 1$. If $s \geq s_0 \sim \log N$, set $\theta_{u,s}(\Delta_s t)=\sqrt{u}\eta_{s+1}^{1/2}\|\Delta_s h\|_{(2u\eta_{s+1})}$.

\begin{Theorem} \label{thm:uncond-est}
There exist absolute constants $c_1$ and $c_2$ for which the following holds. If $\mu$ is an isotropic, unconditional log-concave measure, $H_T=\{\inr{t,\cdot} : t \in T\}$ and $(T_s)_{s \geq 0}$ is an admissible sequence of $T$, then for every $u \geq c_1$,
$$
\theta_{u,s}(\inr{\Delta_s t,\cdot}) \leq c_2u \left(2^s\|\Delta_s t\|_{\ell_\infty^n} + 2^{s/2}\|\Delta_s t\|_{\ell_2^n}\right).
$$
\end{Theorem}

\proof
Let $T \subset \R^n$, and identify it with the class of linear functionals $H_T=\{\inr{t,\cdot} : t \in T\}$ on $(\R^n,\mu)$. By Borell's inequality \cite{Bor}, the $\psi_1$ and $L_2$ norms are $c_1$-equivalent on $\R^n$, where $c_1$ is an absolute constant, and since $\mu$ is isotropic, then $\|\inr{t,\cdot}\|_{L_2} =\|t\|_{\ell_2^n}$. Moreover, there is an absolute constant $c_2$ such that for every $p \geq q$ and $t \in \R^n$,
$$
\|\inr{t,\cdot}\|_{L_p} \leq c_2\frac{p}{q}\|\inr{t,\cdot}\|_{L_q}.
$$
Hence, for every $t \in \R^n$ and every $r \geq 1$,
\begin{align*}
\|\inr{t,\cdot}\|_{(r q)} & \leq \sup_{q \leq \ell \leq r q} \frac{\|\inr{t,\cdot}\|_{L_\ell}}{\sqrt{\ell}}+ \|\inr{t,\cdot}\|_{(q)}
\leq c_2\sup_{q \leq \ell \leq r q}\frac{\ell}{q}\frac{\|\inr{t,\cdot}\|_{L_q}}{\sqrt{\ell}}+\|\inr{t,\cdot}\|_{(q)}
\\
& \leq  (c_2\sqrt{r} +1)\|\inr{t,\cdot}\|_{(q)}.
\end{align*}
Also, for any $t \in \R^n$ and any $p \geq n$, $\|\inr{t,\cdot}\|_{(p)} \leq 2c_2 \sqrt{\frac{p}{n}} \|\inr{t,\cdot}\|_{(n)}$. Therefore,
if $u \geq 1$ and $\eta_s=2^s \leq n$ then
$$
\sqrt{u}\eta_{s+1}^{1/2}\|\inr{\Delta_s t,\cdot}\|_{(2u\eta_{s+1})} \lesssim u 2^{s/2} \|\inr{\Delta_s t,\cdot}\|_{(2^s)},
$$
and if $2^s>n$ then
$$
\sqrt{u}\eta_{s+1}^{1/2}\|\inr{\Delta_s t,\cdot}\|_{(2u\eta_{s+1})} \lesssim \sqrt{u}\frac{2^s}{n}\|\inr{\Delta_s t,\cdot}\|_{(n)}.
$$
Note (\cite{Ball} or \cite{Pa2}, Proposition 3.4) that there is an isotropic convex body $K$ such that for every $t \in \R^n$ and any $1 \leq p \leq n$,
$
\|\inr{t,\cdot}\|_{L_p(\mu)} \leq c_3 \|\inr{t,\cdot}\|_{L_p(K)}
$.
Moreover, since $\mu$ is unconditional, $K$ is also unconditional and using the Bobkov-Nazarov Theorem \cite{BobNaz} we get that
$$ \|\inr{t,\cdot}\|_{L_p(\mu)} \leq c_3 \|\inr{t,\cdot}\|_{L_p(K)}\ls c_{4}  \|\inr{t,\cdot}\|_{L_p(K_{1})}, $$
where $K_{1}$ is an isotropic image of $B_{1}^{n}$.

The moments of every linear functional $\inr{t,\cdot}$ relative to the volume
measure of an isotropic position of $B_1^n$ are well known \cite{GK}: namely, for  $1 \leq p \leq n$,
\begin{equation*}
\|\inr{t,\cdot}\|_{L_p(K_{1})} \sim p\|t\|_\infty +\sqrt{p}\left(\sum_{i=p+1}^n (t_i^2)^* \right)^{1/2}.
\end{equation*}
Combining the two estimates, for $p \leq n$ and any $t \in \R^n$,
\begin{align*}
\|\inr{t,\cdot}\|_{(p)} & \leq  \sup_{q \leq p}
\frac{\|\inr{t,\cdot}\|_{L_q(K_1)}}{\sqrt{q}} \sim  \sup_{q
\leq p} \left(\sqrt{q}\|t\|_\infty +\left(\sum_{i=q+1}^n
(t_i^2)^* \right)^{1/2} \right)
\\
& \leq \sqrt{p}\|t\|_{\ell_\infty^n}+\|t\|_{\ell_2^n}.
\end{align*}
Thus, for $2^s \leq n$,
$$
2^{s/2}\|\inr{\Delta_s t, \cdot}\|_{(2^s)} \leq 2^s \|\Delta_s t\|_{\ell_\infty^n} +2^{s/2}\|\Delta_s t\|_{\ell_2^n},
$$
and if $2^s>n$,
$$
\frac{2^s}{\sqrt{n}}\|\inr{\Delta_s t,\cdot}\|_{(n)} \leq 2^s\|\Delta_s t\|_{\ell_\infty^n}.
$$
\endproof
Note that for an almost optimal admissible sequence,
$$
\sum_{s \geq 0} 2^s\|\Delta_s t\|_{\ell_\infty^n} +2^{s/2}\|\Delta_s t\|_{\ell_2^n} \lesssim \gamma_1(T,\ell_\infty^n) + \gamma_2(T,\ell_2^n).
$$
 It turns out that $\gamma_1(T,\ell_\infty^n) + \gamma_2(T,\ell_2^n)$ can be completely characterized by the following beautiful result due to Talagrand \cite{Tal-AJM,Tal-book}.
\begin{Theorem} \label{Thm:equiv-gamma1-and-gamma-2}
There exist absolute constants $c$ and $C$ for which the following holds. Let $(y_i)_{i=1}^n$ be independent, standard exponential variables. Then, for every $T \subset \R^n$,
$$
c \E \sup_{t \in T} \sum_{i=1}^n y_i t_i \leq \gamma_1(T,\ell_\infty) + \gamma_2(T,\ell_2) \leq C \E \sup_{t \in T} \sum_{i=1}^n y_i t_i.
$$
\end{Theorem}

Recall that if $(y_i)_{i=1}^n$ are standard exponential random variables and $T \subset \R^n$, then we denote $E(T)=\E \sup_{t \in T} \sum_{i=1}^n y_i t_i$ and $d_2(T)=\sup_{t \in T} \|t\|_2$.

Combining the estimates above, it follows that on $\Omega_{1,u} \cap \Omega_{2,u} \cap \Omega_{3,u}$, $P_\sigma T$ satisfies Definition \ref{def:decomposition} with $\theta_s=2^s\|\cdot\|_{\ell_\infty^n}+2^{s/2}\|\cdot\|_{\ell_2^n}$ for $s \geq s_0$ and $\theta_s=0$ otherwise, $\gamma \lesssim E(T)$, $\phi \sim \phi_1$, $\|(\inr{X_i,t})_{i=1}^N\| = \|\inr{t,\cdot} \|_{\psi_1} \sim \|t\|_{\ell_2^n}$ and $\alpha \sim u$. Therefore, $B_4 \lesssim \gamma_2(T,\ell_2) \lesssim E(T)$.
\begin{Theorem} \label{thm:Bernoulli-unconditional}
There exist absolute constants $c_1$, $c_2$, $c_3$ and $c_4$ for which the following holds. For every $u \geq c_1$, With $\mu^N$-probability at least $1-2\exp(-c_2u\log N)$, the set $V=P_\sigma T$ satisfies that
$$
\sup_{v \in V} \left|\sum_{i=1}^N \eps_i v_i^2 \right| \leq c_3ru^{2} \left(d_2(T)\sqrt{N} E(T) + (E(T))^2 \right)
$$
with probability at least $1-2\exp(-c_4 r^2)$ with respect to the Bernoulli random variables.
\end{Theorem}

\subsection{Proofs of Theorems A and B} \label{sec:proofsA-B}
The final step we need for the proofs of Theorem A and Theorem B is a version of the Gin\`{e}-Zinn symmetrization Theorem (see, e.g. \cite{GZ84,VW}), which enables one to pass from the Bernoulli process indexed by random coordinate projections of a class of functions, to the empirical process indexed by the class.

\begin{Theorem} \label{thm:GZ-sym}
Let $F$ be a class of functions and for every $x>0$, set $\beta_N(x)=\inf_{f \in F} Pr(|\sum_{i=1}^N f(X_i) - \E f| > x/2)$. Then
$$
\beta_N(x) Pr_{X} \left(\sup_{f \in F} \left| \sum_{i =1}^N f(X_i) - \E f \right| > x \right) \leq 2Pr_{X \otimes \eps} \left(\sup_{f \in F} \left|\sum_{i=1}^N \eps_i f(X_i) \right| >x/4 \right).
$$
\end{Theorem}
To apply Theorem \ref{thm:GZ-sym}, one has to identify the right value $x$ for which $\beta_N(x) \geq 1/2$. In our case, $F=H^2$, and thus one has to show that if $x$ is large enough, then $\sup_{h \in H} Pr(|\sum_{i=1}^N h^2(X_i) - \E h^2| >x/2) \leq 1/2$.

\begin{Lemma} \label{lemma:check-Gine-Zinn}
Let $H$ be a class of functions which is bounded in $L_q$ and consider the empirical process indexed by $F=\{h^2 : h \in H\}$. If $q \geq 4$ and $x \gtrsim d_{L_q}^2 \sqrt{N}$ then $\beta_N(x) \geq 1/2$ and the same holds if $2<q< 4$ and $x\gtrsim_q d_{L_q}^2 N^{2/q}$.
\end{Lemma}

\proof
The first part of the claim follows from an application of Chebyshev's inequality, and is omitted. For the second part, fix $r>0$, set $V=(h^2(X_i))_{i=1}^N$, and since $Pr(|h^2(X)| \geq d_{L_q}^2 (rN/i)^{2/q}) \leq i/(rN)$ then
\begin{equation} \label{eq:lemma-520}
Pr \left( V_i^* \geq d_{L_q}^2(r N/i)^{2/q}\right) \lesssim \binom{N}{i}\cdot(i/r N)^i \leq \exp(-i\log(er))=(er)^{-i}.
\end{equation}
Moreover, for $r_0 \sim c_2$, $Pr(\exists i : |h^2(X_i)| \geq r_0d_{L_q}^2N^{2/q}) \leq 1/10$. Hence, a truncation argument shows that without loss of generality we may assume that $\|h^2\|_{L_\infty} \leq r_0d_{L_q}^2N^{2/q}$. Applying the $L_\infty$ estimate for the largest two coordinates of $V$ and \eqref{eq:lemma-520} for the rest, it follows that
$$
Pr( \|V\|_{\ell_2^N} \geq c_q(r_0+r)d_{L_q}^2 N^{2/q}) \lesssim \sum_{i=3}^N r^{-i} \lesssim r^{-2}.
$$
Hence, under the truncation assumption,
\begin{align*}
\E |\sum_{i=1}^N h^2(X_i) - \E h^2 | \lesssim & \E_X \E_\eps |\sum_{i=1}^N \eps_i h^2(X_i)| \lesssim \E_X (\sum_{i=1}^N h^4(X_i))^{1/2}=\E \|V\|_{\ell_2^N} \\
\lesssim_q & r_0d_{L_q}^2 N^{2/q},
\end{align*}
showing that it suffices to take $x \sim_q d_{L_q}^2 N^{2/q}$ as claimed.
\endproof

Since $x/N$ is well within our range, one may complete the proofs of Theorem A and Theorem B.

\noindent{\bf Proof of Theorem A.}
For $q>4$, let $\rho_{r,u} \sim_{\kappa_1,\kappa_2,q} ru (\sqrt{n/N}+n/N)$ for $u \gtrsim_{\kappa_1,\kappa_2,q} c_1$, and $r \geq c_2$. If $2<q<4$ set $\rho_{r,u} \sim_{\kappa_1,\kappa_2,q} ru (n/N)^{1-2/q}$ for $u \gtrsim_{\kappa_1,\kappa_2,q} \log(eN/n)$ and $r \geq c_3$. Then,
\begin{align*}
& Pr_X \left(\sup_{t \in S^{n-1}} |\frac{1}{N}\sum_{i=1}^N \inr{t,X_i}^2 -\E \inr{t,X}^2 | \geq \rho_{r,u}\right)
\\
\leq & 4\E_X Pr_\eps \left(|\frac{1}{N}\sum_{i=1}^N \eps_i\inr{t,X_i}^2| \geq \rho_{r,u}/4\right)
\\
\lesssim & Pr_X ((\Omega_{1,u}\cap \Omega_{2,u} \cap \Omega_{3,u})^c) +  Pr_\eps  \left(|\frac{1}{N}\sum_{i=1}^N \eps_i\inr{t,X_i}^2| >\rho_{r,u}/4 \Big{|} \Omega_{1,u} \cap \Omega_{2,u} \cap \Omega_{3,u} \right)
\\
\lesssim & \exp(-c_4n).
\end{align*}
\endproof

\noindent {\bf Proof of the quantitative Bai-Yin Theorem.}

To prove the quantitative version of the Bai-Yin Theorem one has to combine Theorem A with a conditioning argument. Consider the vector $X=(\xi_1,...,\xi_n)$ with $\xi \in L_q$ for some $q>4$, and let $\nu$ be the measure on $\R^n$ given by $\nu=X | cn^{1/p}B_p^n$; that is, $\nu$ is given by the conditioning of $X$ to the unconditional body $cn^{1/p}B_p^n$ for a suitable choice of $c$ and $p$. Clearly, $\nu$ is unconditional and satisfies the $p$-small diameter $L_q$ moment assumption, and thus, falls within the realm of Theorem A.  Therefore, if the event ${\cal A}=\{\max_{i \leq N} \|X_i\|_{\ell_p^n} \leq cn^{1/p}\}$ has high enough probability, the quantitative version of the Bai-Yin Theorem follows from Theorem A, because for every event ${\cal B}$,
$$
Pr( (X_i)_{i=1}^N \in {\cal B}) \leq Pr( (X_i)_{i=1}^N \in {\cal B} | X_1,...,X_N \in cn^{1/p}B_p^n) Pr({\cal A}) +Pr({\cal A}^c).
$$
Hence, the final step in the proof of our version of the Bai-Yin Theorem is to show that if $\xi \in L_q$ for $q>2$, there is some $p>2$ for which ${\cal A}$ has a large measure.

Recall that for every $v \in \R^n$, $\|v\|_{\ell_{p,\infty}^n}=\max_{i \leq n} v_k^*/k^{1/p}$, and since $\ell_r^n \subset \ell_{p,\infty}^n \subset \ell_p^n$ for every $r<p$, it suffices to show that $\max_{i \leq N} \|X\|_{\ell_{p,\infty}^n} \lesssim n^{1/p}$ for some $p>2$ with high enough probability.
\begin{Lemma} \label{lemma-independent-small-diameter}
For every $q>4$ and $2<p<q$, there exist constants $c_1$ and $c_2$ that depend on $q$ and $p$ for which the following holds.
If $\xi \in L_q$, $X=(\xi_1,...,\xi_n)$ and $X_1,...,X_N$ are independent copies of $X$, then
$$
Pr(\max_{1 \leq i \leq N} \|X_i\|_{\ell_{p,\infty}^n} \geq c_1\|\xi\|_{L_q} n^{1/p}) \leq \frac{c_2N}{n^{\frac{q}{p}-1}}.
$$
\end{Lemma}

\proof
If $A=\|\xi\|_{L_q}$ then $Pr(|\xi| \geq At) \leq t^{-q}$, and for every $1 \leq k \leq n$, $Pr(\xi_k^* \geq t) \leq \binom{n}{k}(Pr(|\xi| \geq t))^k$. Therefore, if $p<q$ and $y>e$ then
\begin{align*}
& Pr(\xi_k^* \geq A(ny/k)^{1/p}) \leq \exp(k\log(en/k)-k(q/p)\log(ny/k))
\\
\leq & \exp(-k(\frac{q}{p}-1)\log(ny/k)).
\end{align*}
Using this estimate for every $k=2^j$ and summing the probabilities, it follows that for every $q$ and $p$ there is a constant $c_{q,p}$ for which $\|X\|_{\ell_{p,\infty}^n} \lesssim n^{1/p}$ with probability at least $1-c_{q,p}n^{1-q/p}$, and in particular, $Pr(\max_{i \leq N} \|X_i\|_{\ell_{p,\infty}^n} \geq cn^{1/p}) \leq c_{q,p}N/n^{(q/p)-1}$, as claimed.
\endproof
Combining Lemma \ref{lemma-independent-small-diameter} with Theorem A concludes the proof of the quantitative Bai-Yin Theorem.
\endproof

\noindent{\bf Proof of Theorem B.}
If $r \sim u$, with probability at least $1-2\exp(-c_3u^2)$ with respect to the Bernoulli random variables,
$$
\sup_{v \in P_\sigma T} \left|\frac{1}{N}\sum_{i=1}^N \eps_i v_i^2 \right| \lesssim_{u^3} d_2(T)\frac{E(T)}{\sqrt{N}}+\frac{E^2(T)}{N}.
$$
Since $d_2(T)E(T)/\sqrt{N}$ is a ``legal" choice in the Gin\'{e}-Zinn symmetrization theorem, the proof is concluded.
\endproof

\footnotesize {

\end{document}